\documentclass{article}
\usepackage[T2A]{fontenc}
\usepackage[cp1251]{inputenc}
\usepackage[english,russian]{babel}
\usepackage[tbtags]{amsmath}
\usepackage{amsfonts,amssymb,mathrsfs,amscd}

\let\No=\textnumero

\usepackage[hyper]{izv2e}
\JournalName{}

\numberwithin{equation}{section}
\theoremstyle{plain}
\newtheorem{theorem}{Теорема}
\newtheorem{lemma}{Лемма}[section]
\newtheorem{propos}{Предложение}
\newtheorem{corollary}{Следствие}
\theoremstyle{definition}
\newtheorem{definition}{Определение}
\newtheorem{proof}{Доказательство}
\newtheorem{remark}{Замечание}

\begin{document}
	
\title{ Models of representations for classical series of  Lie algebras }
\author[D.\,V.~Artamonov]{D.\,V.~Artamonov}
\address{Lomonosov Moscow State University}
\email{artamonov.dmitri@gmail.com}

\date{19.03.2008}
\udk{517.588}

\maketitle

\begin{fulltext}
	
	\begin{abstract}
A model of representations of a Lie algebra is a representation which a direct sum of all irreducible finite dimensional representations taken with  the multiplicity  $1$. 
In the paper an explicit construction of a model of representations for all  classical series of  simple  Lie algebras is given. The construction does not differ much for different series. The space of the model is constructed as a space of polynomial solutions of a system of partial differential equations. The equations in this system are constructed form relations between minors of matrices taken from the corresponding Lie group.  This system has a simplification which is very close to the GKZ system, that is satisfied by
$A$-hypergeometric functions.
		
	Bibliography: 32 items.
	\end{abstract}
	
	
	\markright{А-ГКЗ модели представлений}
	
	
	\section{Введение}
	
	A model of representations of a Lie algebra is a representation which is a direct sum of all it's finite dimensional irreducible representations taken with multiplicity    $1$.

	Thus one can think about  the classical Weyl construction as a model of representations of the algebra  $\mathfrak{gl}_n$  \cite{w}.  Formally the Weyl construction is an explicit embedding of a certain  irreducible representation into a  tensor power of a standard representation of  $\mathfrak{gl}_n$. But when one take a direct sum of these embeddings one obtains a model of representations in the sense of the present paper.  There exist analogs of the Weyl construction for other classical Lie algebras  \cite{fh}, and even for some exceptional Lie algebras \cite{hu}.
	In physical literature the models based on the language of creation and annihilation operators  are used. Such an approach in the case of the series   $A$ is used in \cite{bb}  and further papers of these authors. But when one tries to generalize this language to the the case of the series $C$ one faces difficulties. Thus usually small dimensions are considered   \cite{h}, \cite{h1}, \cite{h2}, \cite{h3}.
	Also the Zhelobenko's construction is a model  \cite{zh}.

	Actually all these three constructions are similar. The construction of the present paper is also similar to them.
	
	Let us continue  the discussion of known models. There exist numerous models of combinatorial  nature, it seems to be impossible to give short  their  review.  Let us return to the mentioned above papers by  Biedenharn and coauthors  \cite{bb}-\cite{h3}.
	In paper  \cite{bf} (see also  \cite{bf1}) there was constructed a model of representations for the algebra    $\mathfrak{sl}_3$.   These two papers were the closing papers  in a  big series of papers where  the authors tried to  calculate explicitly the Clebsh-Gordan, Rachah coefficients that describe a splitting  of tensor products of representations into irreducible summands. Flath joint to this activity on it's final stage and although he dealt with  the classical objects he named the explicit classical representation theory the  "mathematical golden mine"  \cite{bf1}.
	Inspired by  \cite{bf},  Gelfand and Kapranov wrote a  paper   \cite{gz}, where the notion of a model of representations was introduced.
	In this paper some models of geometric nature for all classical  Lie algebras over  $\mathbb{C}$ were constructed. These written with an inspiration papers devoted to the classical representation theory were the starting point of the present paper.

It is necessary to mention the existence  of other numerous models of geometric nature. Thus the exist models realized as subspaces  in the space of functions on a homogeneous space, on a  	   $HV$-variety	 \cite{VP}

	The present paper can be considered as a continuation of the paper   \cite{a4} and as a generalization of it's result to other classical Lie algebras.
	In  \cite{a4}  for the algebra  $\mathfrak{gl}_n$  we consider the function of independent variables  $A_X$, $X\subset\{1,...,n\}$, anti-symmetric on  $X$. We construct a model of representations formed by  polynomials in these variables, that satisfy some system  of partial differential equations, called the antisymmetrized system of Gelfand-Kapranov-Zelevinsky  (A-GKZ for short, see \cite{a5}).  This model is called the A-GKZ model.  Note that this system of equations  is close to the hypergeometric system on the space   $\Lambda^k\mathbb{C}^N$,constructed in  \cite{GG}. Hence  it is also can be called a system of hypergeometric type. 
	
Thus one obtains a model of representations whose space is a space of polynomial solutions of a system of hypergeometric type. The existence of such a model can be considered as explanation of the fact that in explicit calculations in the representation theory very often hypergeometric functions and constants appear (see for example  \cite{vil}).

	 Moreover in this model one naturally constructs a base in each representation. An advantage of this model is the fact that one has both an explicit base and an explicit formulas for the scalar product. This makes possible to do some non-trivial calculations. Thus in \cite{sm},  using the A-GKZ model in the case  $n=3$ explicit simple formulas for an arbitrary Clebsh-Gordan coefficient are obtained.   And in the paper  \cite{a6j} using the same model an explicit formula for a $6j$-symbol for the algebra $\mathfrak{gl}_3$ is obtained.

	Also in  \cite{a4} using   $A$-hypergeometric functions a base in Zhelobenko's realization was obtained (in it's construction the A-GKZ system plays a crucial role). In this base one manages to write explicitly the action of generators of the algebra  $\mathfrak{gl}_n$. This construction is called the GKZ base for the Zhelobenko's model.  In   \cite{a4} this base is an important tool in establishing a relation  between the constructed base in the  A-GKZ model and the Gelfand-Tsetlin base. 
	
	There exist analogs of the GKZ base for other Lie algebras of small dimensions \cite{a3}, \cite{a6}.

	In the present paper in Section  \ref{systemy}  we construct an analogue of the GKZ and A-GKZ systems for the Lie algebras of the series   $B$, $C$, $D$.  This leads to an A-GKZ model for these algebras (see Section  \ref{modela}).  It is remarkable that the constructions for different series are essentially the same.
	In the A-GKZ model one naturally constructs a base.

	Also for the series $B$, $C$, $D$  we construct a GKZ base in the Zhelobenko's model and establish formulas for the action of generators in this base  (see  Section \ref{zhlb}). Using this result in Section   \ref{monomy} we construct other basis in the Zhelobenko's realization.

	Also we investigate a relation between the constructed base in the A-GKZ model and the Gelfand-Tsetlin base. 
	
	Firstly we suggest a new point of view to the notion of the Gelfand-Tsetlin diagram  (see Section  \ref{dgc},  Definition  \ref{d1}).
	There exists a bijective correspondence between the objects introduced in Definition  \ref{d1} and the traditional Gelfand-Tsetlin diagrams.  In our approach the formulas for the action of generators and so on are written more natural (see Section \ref{gkzt}).

	The Gelfand-Tsetlin diagrams in our sense  are indexing base vectors in irreducible representations. Of course there are numerous constructions of sets that solve the same problem.  Usually they are constructed as sets of integer points in some polytops (the string polytops of Berenstein-Zelevinsky-Littelmann \cite{abz}, the polytops of Vinberg-Littelmann-Feigin-Fourier  \cite{ffl}).
	The objects introduced in   \ref{d1}  can also be obtained as integer points in some polytop. But we use the term the Gelfand-Tsetlin diagram since one can easily reconstruct the highest weights for a chain of subalgebras  that appear  in the standard procedure of construction of the Gelfand-Tsetlin base (diagrams for classical series in the traditional  sense can be found in\cite{m}).

	Secondly we prove that the transition matrix between the base of the A-GKZ model and the Gelfand-Tsetlin base is triangular. We show that the Gelfand-Tsetlin base is nothing but an orthogonalization of the A-GKZ base.
	

	The structure of the paper is the following. Section  \ref{oob}, \ref{km}, \ref{seriaa}  are introductory. In  \ref{oob}  the  basic notions are introduce. In \ref{km} the Zhelobenko's model is discussed, here we obtain an explicit description of the Zhelobenko's model for the series   $B$, $C$, $D$,  supplementing the results of the book  \cite{zh}.  The results in this Section are formally new but they can be obtained form results  of the book  \cite{zh}  by simple calculations.

	In  \ref{seriaa} we present the results for the series  $A$, which are modifications of the results from  \cite{a4}, there are no essential new results in this Section.

	The main result of the paper can be found in Sections  \ref{systemy}, \ref{dgc}, \ref{mdl}. In   \ref{systemy}  the A-GKZ system for the series $B$, $C$, $D$ is introduced, in Section  \ref{dgc}  a new definition of a Gelfand-Tsetlin diagram for this series is given,  in Section \ref{mdl}  the A-GKZ model and the GKZ base for the Zhelobenko's model are constructed.
	
In  Section \ref{gttl}  a relation of the constructed base of the A-GKZ model and the Gelfand-Tsetlin base  is discussed.

	Finally let us note the paper 	\cite{obl}. The results  of this paper are not directly related to the results of the present paper, but the ideologically  they are quite close. Here    {\it  an arbitrary   }  GKZ function is interpreted as a matrix element in a representation of a special Lie algebra.
	
	\section{The basic objects}
	\label{oob}
	
	In this Section we give an definition of an important class of functions that plays a crucial role in the present paper. We present a system of equations satisfied by these functions. Also we present the Lie algebras used in the present paper. 
	
	\subsection{ A $\Gamma$-series}
	\label{r2}
	
	A detailed information about a $\Gamma$-series can be found  in \cite{GG}.
	
	Let $\mathcal{B}\subset \mathbb{Z}^N$ be a lattice, let $\gamma\in \mathbb{Z}^N$ be a fixed vector. Define a {\it  hypergeometric 
		$\Gamma$-series }  in variables  $z_1,...,z_N$ by the formula
	
	\begin{equation}
		\label{gmr}
		\mathcal{F}_{\gamma}(z,\mathcal{B})=\sum_{b\in
			\mathcal{B}}\frac{z^{b+\gamma}}{\Gamma(b+\gamma+1)},
	\end{equation}
	where  $z=(z_1,...,z_N)$.  The following  multi-index notation  is used
	
	$$
	z^{b+\gamma}:=\prod_{i=1}^N
	z_i^{b_i+\gamma_i},\,\,\,\Gamma(b+\gamma+1):=\prod_{i=1}^N\Gamma(b_i+\gamma_i+1).
	$$
	
	Note that  if at least one of the components of the vector  $b+\gamma$ is negative integer, then the corresponding summand in     \eqref{gmr} vanishes. Due to this fact in the $\Gamma$-series, considered in the present paper,    there are finitely many non-zero terms.  For simplicity we shall write factorials instead of $\Gamma$-functions.

	An $A$-hypergeometric function satisfies a system of partial differential equations called the Gelfand-Karpanov-Zelevinsky (shortly GKZ) system, which consists of equations of two types.
	
	{\bf 1.} Let  $a=(a_1,...,a_N)$ be a vector orthogonal to the lattice  $\mathcal{B}$, then
	
	\begin{equation}
		\label{e1}
		a_1z_1\frac{\partial}{\partial z_1}\mathcal{F}_{\gamma}+...+a_Nz_N\frac{\partial}{\partial z_N}\mathcal{F}_{\gamma}=(a_1\gamma_1+...+a_N\gamma_N)\mathcal{F}_{\gamma},
	\end{equation}
	it is enough to consider only the base vectors of the lattice orthogonal to  $\mathcal{B}$.

	{\bf 2.} Let  $b\in \mathcal{B}$ Рё $b=b_+-b_-$, where all coordinates of the vectors  $b_+$, $b_-$ are non-negative.  Chose in these vectors non-zero elements  
	$b_+=(...b_{i_1},....,b_{i_k}...)$,  $b_-=(...b_{j_1},....,b_{j_l}...)$. Then
	
	\begin{equation}
		\label{e2} (\frac{\partial }{\partial
			z_{i_1}})^{b_{i_1}}...(\frac{\partial}{\partial z_{i_k}})^{b_{i_k}}
		\mathcal{F}_{\gamma}=(\frac{\partial }{\partial
			z_{j_1}})^{b_{j_1}}...(\frac{\partial }{\partial z_{j_l}})^{b_{j_l}} \mathcal{F}_{\gamma}.
	\end{equation}
	
It is enough  to consider only a finite collection fo vectors   $b\in \mathcal{B}$ \footnote{In the paper  \cite{GG}  it is  wrongly stated that it is enough to consider the equation corresponding  to base vectors of the lattice  $B$.}.  The system of partial differential  equations can be identified with an ideal in  the  ring of differential operators, generated by the operators defining the equations of the system.  Let us present the way the construct  explicitly the ideal corresponding to the system  \eqref{e2} .  To do it let us firstly present the properties of the correspondence  $b\in B\mapsto \mathcal{O}_b$. One has:

\begin{enumerate}
	\item $kb\mapsto (\mathcal{O}_b)^k$, $k\in\mathbb{Z}_{\geq 0}$, $\,\,\,\,$ $-b \mapsto -\mathcal{O}_b$,
	\item   Let   $b=b_{+}-b_{-}$, $c=c_{+}-c_{-}$,  where  $c_{\pm}$ have only non-negative coordinates. Assume that  the presentation of  $b+c$  as a difference of vectors with non-negative coordinates is $(b_{+}+c_{+})-(b_{-}+c_{-})$ (that is in this equality in each coordinate no reduction takes place). Then
	
	$$
	b+c\mapsto \mathcal{O}_{b+c}=(\frac{\partial}{\partial z})^{c_+}\mathcal{O}_{b}+(\frac{\partial}{\partial z})^{b_+}\mathcal{O}_{c}
	$$

	\item  Let   $b=(b_{+}+u)-(b_{-}+v)$, $c=(c_{+}+v)-(c_{-}+u)$,  where  $c_{pm},u,v$ have only non-negative coordinates. Assume that  the presentation of  $b+c$  as a difference of vectors with non-negative coordinates is  $(b_{+}+c_{+})-(b_{-}+c_{-})$ (that is in this equality in each coordinate no reduction takes place). Then
	
	$$
	(\frac{\partial}{\partial z})^{u+v} \mathcal{O}_{b+c}=(\frac{\partial}{\partial z})^{c_+}\mathcal{O}_{b}+(\frac{\partial}{\partial z})^{b_+}\mathcal{O}_{c}
	$$
\end{enumerate}

\begin{definition}\label{porozd}
Let us be given a collection of differential operators with constant coefficients. A system    {\it  generated  }  by this collection is a system defined by the following collection of differential operators. Firstly we take all operators form the initial collection. Secondly all operators the  belong to the ideal generated by these operators. Thirdly all operators  that are obtained from the operators belonging to the obtained ideal by  division (if it is possible) by a differential monomial.
\end{definition}

The considerations above show that the  GKZ system is   {\it  generated  } by operators   corresponding to base vectors of the lattice.

	\subsection{ Algebras $\mathfrak{o}_{2n+1}$,  $\mathfrak{o}_{2n}$, $\mathfrak{sp}_{2n}$}
	\label{algl}
	
	The Lie algebras  $\mathfrak{o}_{2n}$, $\mathfrak{sp}_{2n}$ are considered as subalgebras in the Lie algebra of all  $2n\times 2n$ matrices,   whose rows and columns are indexed by   $i,j=-n,...,-1,1,...,n$, and the algebra   $\mathfrak{o}_{2n+1}$   is a subalgebra  in the Lie algebra of all $(2n+1)\times (2n+1)$ matrices, whose rows and columns are indexed by    $i,j=-n,...,-1,0,1,...,n$.

	The algebras $\mathfrak{o}_{2n+1}$ and   $\mathfrak{o}_{2n}$
	are generated by matrices
	
	\begin{equation}\label{fbd}F_{i,j}=E_{i,j}-E_{-j,-i},\end{equation}
	
	where $i,j=-n,...,-1,0,1,...,n$ in the case $\mathfrak{o}_{2n+1}$ and  $i,j=-n,...,-1,1,...,n$   in the case $\mathfrak{o}_{2n}$.
	
	The algebra $\mathfrak{sp}_{2n} $  is generated by the matrices
	\begin{equation}\label{fd}F_{i,j}=E_{i,j}-sign(i)sign(j)E_{-j,-i},\end{equation}
	where $i,j=-n,...,-1,1,...,n$.

	Denote the Lie algebras  $\mathfrak{o}_{2n+1}$,  $\mathfrak{o}_{2n}$ or  $\mathfrak{sp}_{2n}$ just as $g_n$.

	The chosen realization is also a root realization. The elements    $F_{i,j}$, $i<j$ correspond to positive roots;  $F_{i,j}$, $i>j$ correspond to negative roots; $F_{i,i}$ generate the Cartan subalgebra.
	
An analogous choice of the Cartan subalgebra and the root elements is used for 	  $\mathfrak{gl}_{m}$.
	
	To be able to talk about a Gelfand-Tsetlin base one needs to fix a chain of subalgebras. 
	A subalgebra $g_{n-k}\subset g_{n}$ is defined as a span  $<F_{i,j}>_{i,j\neq \pm 1,...,\pm (k-1)}$.

	\section{ The Zhelobenko's model}
	\label{km}
	
	In this Section we present a model of representations realized in the space of functions  on the corresponding Lie group. Zhelobenko (see  \cite{zh})  has proved the Theorem  \ref{tzh0}  which describes the space of this model in the case of the series  $A$. We formulate and prove the Theorem   \ref{tzh},  which is a generalization of the Zhelobenko's result to the cases of other series. 
	
	In these Theorems the space of a representation is described as a space of solution of  a system of partial differential equations. A more explicit description is given in the Theorem \ref{ost1}.

	\subsection{Functions on a group}

	Consider a space of function on a group
	$G$. Onto a function
	$f(g)$,  $g\in G$ an element $X\in G_{}$ acts by right shifts  
	\begin{equation}
		\label{xf}
		(Xf)(g)=f(gX).
	\end{equation}
	
	Thus the space  $Fun$ of all function of a group    $G$ is a representation of   $G$ and thus of the Lie algebra    $Lie G$.
	
	Let  $G=Sp_{2n}$, $SO_{2n+1}$, $SO_{2n+1}$  and let
	$a_{i}^{j}$, $i,j=-n,...,n$\footnote{In the cases $G=Sp_{2n}$, $O_{2n}$ one takes $i,j=-n,...\hat{0}...,n$, and in the case  $G=SO_{2n+1}$  one takes $i,j=-n,...,0,...,n$} 
	be a function of a matrix element which is a function on the group $G$. Here $j$ is a row index, and  $i$ is a column index.

	Also put
	
	\begin{equation}
		\label{dete}
		a_{i_1,...,i_k}:=det(a_i^j)_{i=i_1,...,i_k}^{j \in\text{the first  $k$ rows}}.
	\end{equation}

	That is one takes a determinant of a submatrix in a matrix $(a_i^j)$,
	formed by the first  $k$ rows and columns  
	$i_1,...,i_k$

	Also in the case of the series   $D$ one puts
	\begin{equation}
		\bar{a}_{i_1,...,i_n}:=det(a_i^j)_{i=i_1,...,i_n}^{j=-n,...,-2,1}.
	\end{equation}

	One has
	
	\begin{equation}
		(a_{-n,...,-2,-1})^{-1}= \bar{a}_{-n,...,-2,1}.
	\end{equation}
	
	Using  formulas \eqref{xf}, \eqref{dete} one can obtain formulas for the action of generators of the algebra onto these determinants.  To do it let us formally \footnote{The expression  \eqref{edet1} should be understood symbolically since it does not define an action of    $E_{i,j}$, because it does not respect the relations between determinant in the cases of groups   $Sp$,  $SO$.} introduces an action of an operator
	$E_{i,j}$ onto a determinant:
	
	\begin{equation}
		\label{edet1}
		E_{i,j}a_{i_1,...,i_k}=a_{\{i_1,...,i_k\}\mid_{j\mapsto i}},
	\end{equation}

	where  $.\mid_{j\mapsto i}$ is a substitution of   $j$ instead of
	$i$, and in the case when  $j\notin \{i_1,...,i_k\}$, then one obtains zero.
	Now for the Lie algebra of the series  $B, C, D$ one defines an action of  $F_{i,j}$  through the action of  $E_{i,j}$.

	Now let us give an explicit formulas for the highest vector of a given highest weight  (see  \cite{bb}, 
	\cite{zh}). The vector
	
	\begin{align}\begin{split}
			\label{stv0} &v_0=\prod_{k=-n}^{-2} (a_{-n,...,-k})^{m_{-k}-m_{-k+1}}
			a_{-n,...,-2,-1}^{m_{-1}} \text{  for the series $B$, $C$, $D$ in the case
				$m_{-1}\geq 0$},\\
			&v_0=\prod_{k=-n}^{-2} (a_{-n,...,-k})^{m_{-k}-m_{-k+1}}
			\bar{a}_{-n,...,-2,1}^{-m_{-1}} \text{ for the series  $D$ in the case $m_{-1}<
				0$}
	\end{split}\end{align}

	is a highest vector for  $g_{n}$ with the highest weight
	$[m_{-n},...,m_{-1}]$ for all the series  $B$, $C$, $D$.
	Note that in the case of the integer highest weight this is a polynomial function. In the case  $B$, $D$
	when the highest weight is half-integer the last factor has a fractional exponent.  Thus one obtains a multi-valued function. It becomes single-valued when one passes toe the group   $Spin$.

	\begin{definition}
		A direct sum of subrepresentations in  $Fun$,  generated by highest vectors \eqref{stv0}  is called the Zhelobenko's model. The space of this model is denoted by     $Zh$.
	\end{definition}

	\subsection{Relations between determinant}
	
	\subsubsection{ The Plucker relations}
	
	For all series between the minors $a_{i_1,...,i_k}$ of an  $m\times m$  matrix one has  the Plucker relations 
	
	\begin{equation}
		r_s=\sum_{t=1}^{k+1} (-1)^t a_{i_1,...,i_{k-1},j_t} a_{j_1,...,j_{t-1}, j_{t+1},...,j_{k+1}}=0,
	\end{equation}
	
where $s$ is an index numerating these relations.

	These relations are sufficient conditions that guarantee that the collection of numbers  $a_{i_1,...,i_k}$ is a collection of minors of type \eqref{dete} of some matrix.
	That is one has.
	
	\begin{lemma}[see  \cite{rel}]
		\label{lj0}
		In the case of the series 	  $A$   the ideal  $I_{\mathfrak{gl}_m}$ of relations between the determinants   $a_X:=a_{i_1,...,i_k}$ is generated by the Plucker relations.
	\end{lemma}

	\subsubsection{The Jacobi relations}
	
	In the case of the series $B$, $C$, $D$  there exist other relations.
	To derive them let us use the fact that between minors of a matrix   $(a_i^j)$ and it's inverse $((a^{-1})_i^j)$ there exist the Jacobi relations \cite{gm}:
	
	\begin{equation}
	\label{soja}
	a_{i_1,...,i_k}^{j_1,...,j_k}=det(a)(-1)^{\sum_{p=1}^k i_p+\sum_{q=1}^k j_q}(a^{-1})^{\widehat{i_1},...,\widehat{i_k}}_{\widehat{j_1},...,\widehat{j_k}}.
\end{equation}

Here when one composes a minor form the right side of the equality on  takes all columns with indices  $\{-n,...,n\}$ except  $i_1,...,i_k$, all row with indices   $\{-n,...,n\}$ except    $j_1,...,j_k$.

	Let us write the corollaries of these relations in the cases   $B$, $C$, $D$. Let  $X$ be an element of the corresponding Lie group.  Then

	\begin{align}
		\begin{split}
			\label{omg}
			& X^t\Omega X=\Omega \Leftrightarrow X^{-1}=\Omega^{-1}X^t\Omega,\,\,\,\,\,  \Omega=(\omega_{i,j}),
			\text{ where }\\
			&\omega_{i,j}=\begin{cases} +1,\,\, i=-j, i<0\\ -1,\,\, i=-j, i>0\\ 0 \text{  otherwise }  \end{cases} \text{ for the series  }C,\,\,\,\,\,\,\,\,\,
			\omega_{i,j}=\begin{cases} +1,\,\, i=-j,\\  0 \text{ otherwise }  \end{cases} \text{ for the series  }B,D
		\end{split}
	\end{align}
	
	Note that  $\Omega^{-1}=-\Omega$  in the case of the series  $C$  and  $\Omega^{-1}=\Omega$  for the series   $B$, $D$.
	
	
	Thus for the matrix elements in the considered case one has 
	\begin{align*}
		&a_{i}^j=a_{-j}^{-i}\text{  for the series  $B,D$},\,\,\,\,\,\,\,a_{i}^j=sing(i)sign(j)a_{-j}^{-i}\text{ for the series  $C$}
	\end{align*}

	Since for all groups one has $det(a)=1$, we obtain:
	
	\begin{align}
		\begin{split}
			\label{rwo0}
			& a_{i_1,...,i_k}^{-n,...,-n+k-1}=(-1)^{i_1+...+i_k}(-1)^{-n+...+(-n+k-1)}(a^{-1})^{\widehat{i_1},...,\widehat{i_k}}_{\widehat{-n},...,\widehat{-n+k-1}}=\\
			&=(-1)^{i_1+...+i_k}(-1)^{\frac{(-2n+k-1)k}{2}}\cdot\begin{cases} (a_{-j}^{-i})^{\widehat{i_1},...,\widehat{i_k}}_{\widehat{-n},...,\widehat{-n+k-1}}    \text{ for the series  $B$, $D$},  \\ (sing(i)sing(j)a_{-j}^{-i})^{\widehat{i_1},...,\widehat{i_k}}_{\widehat{-n},...,\widehat{-n+k-1}}    \text{ for the series  $C$}  \\   \end{cases}=\\&
			=(-1)^{i_1+...+i_k}(-1)^{\frac{(-2n+k-1)k}{2}}\cdot\begin{cases} (a_{i}^{j})_{\widehat{-i_1},...,\widehat{-i_k}}^{\widehat{n},...,\widehat{n-k+1}}    \text{ for the series  $B$, $D$},  \\ (sing(i)sing(j)a_{i}^{j})_{\widehat{-i_1},...,\widehat{-i_k}}^{\widehat{n},...,\widehat{n-k+1}}    \text{ for the series  $C$} \end{cases}
		\end{split}
	\end{align}
	
	Thus we have proved the following statement. 
	
	\begin{lemma}
		\label{lj}
		For the series  $B$, $C$,  $D$ one has
		
		\begin{equation}
			\label{sb}
			a_{i_1,...,i_k}=\pm  a_{\widehat{-i_1},...,\widehat{-i_k}}
		\end{equation}
		For the series $D$  one also has
		
		\begin{equation}
			\label{sd}
			\bar{a}_{i_1,...,i_k}=\pm \bar{a}_{\widehat{-i_1},...,\widehat{-i_k}}.
		\end{equation}
		
		Here
		
		\begin{equation}
			\label{znk}\pm =s\cdot(-1)^{i_1+...+i_k}(-1)^{\frac{(-2n+k-1)k}{2}},
		\end{equation}  where  $s=1$ in the cases  $B$, $D$  and   $s$ equals $-1$ to the power which equals to the number of rows and columns with negative indices in the case of the series  $C$.
		
	\end{lemma}
	
	The statement \eqref{sd} for the series $D$  can be obtained using reordering or rows.

	Introduce a shorter notation for determinants. If  $X\subset \{-n,...,n\}$, then $a_X:=a_{i_1,...,i_k}$. Analogously introduce a notation   $\bar{a}_X$ in the case of the algebra  $\mathfrak{o}_{2n}$ and  $|X|=n$. In these notations  \eqref{sb} one writes
	
	$$
	a_X=\pm a_{\widehat{-X}}.
	$$
	where for  $X={i_1,...,i_k}$ one puts  $-X:=\{-i_1,...,-i_k\}$ and for   $Y=\{j_1,...,j_k\}$  one puts  $\widehat{Y}:=\{-n,...,n\}\setminus Y$.
	
	One has the following statement.
	
	\begin{lemma}
		\label{lms}  Let  $G$ be a group $Sp_{2n}$, $SO_{2n+1}$, $SO_{2n+1}$.  That is a group of matrices that preserve the bilinear form with the matrix $\Omega=(\omega_{i,j})$,  where $\omega_{i,j}$ are defined in  \eqref{omg}.
		Then the ideal of relations between the determinants $a_X$ of type \eqref{dete} is generated by the Plucker and the Jacobi relations.  The same is true in the case of the series  $D$,  where instead of determinants  $a_X$ of the order  $n$  one takes the determinans $\bar{a}_X$.
		
	\end{lemma}
	\begin{proof}
		
		Firstly consider the determinants  $a_X$.
		One has.
		
		\begin{propos}
			\label{prp}
			Let us be given matrices  $O_1$ and $O_2\in G$. If all their minors constructed on columns that belong to an arbitrary subset   $X$ and first consecutive rows  (i.e. minors of type \eqref{dete}) are equal, then $O_1=TO_2$, where $T$ is a low-unitriangular matrix. Analogously  if all their minors constructed on  rows that belong to an arbitrary subset  $X$ and first consecutive columns then $O_1=O_2T$, where  $T$ is an upper-unitriangular matrix.
		\end{propos}
		
	 A close (but not exactly equivalent)  statement can be found in   \cite{rel} (Proposition 14.2).
	\begin{proof}
		Let us prove the first statement, the second can be proved analogously.
		
	The following well-known fact takes place:  a $k$-dimensional subspace   $L=<x^1,...,x^k>\subset\mathbb{C}^N$,  $x^i=(x^i_1,...,x^i_N)$,  is uniquely defined by it's Plucker coordinates   $\{a_{j_1,...,j_k}=det(x^i_j)^{i=1,..,k}_{j=j_1,...,j_k}\}$.
		
	Take matrices 	  $O_1$, $O_2$.  For an arbitrary   $k$ let  $x^1,...,x^k$ be first  $k$ rows of the matrix   $O_1$, and let   $y^1,...,y^k$ be first   $k$ rows of the matrix  $O_2$.  From the statement formulated in the previous paragraph one obtains that   $<x^1,...,x^k>=<y^1,...,y^k>$. Hence   $O_1=TO_2$, where  $T$ is low-triangular.  Since determinants of submatrices of $O_1$  and $O_2$ constructed on the first consecutive rows and columns  are equal, the matrix  $T$ is low-unitriangular .
		
	\end{proof}

		\begin{corollary} 
			\label{sl1}
			Let for minors of $O$ one has 
			$a_{X}=\pm a_{\widehat{-X}}$, where the sign is defined in   \eqref{znk}. Then one has  $\pm\Omega^{-1} O^t\Omega=O^{-1}T$,  
			where $T$ is an upper-unitriangular matrix and $\pm=-$ for the series  $C$ and  $+$ for the series  $B$, $D$.
			  \end{corollary}
		
	\begin{proof}		
		
		The relation $a_{X}=\pm a_{\widehat{-X}}$  was obtained in the formula \eqref{rwo0} as an equality between minors constructed on the first consecutive   {\it columns } for matrices
		$X^{-1}$ and $\pm\Omega^{-1}X^t\Omega$.  Then the Corollary \ref{sl1}  is obtained immediately by application of Proposition  \ref{prp}.

	\end{proof}	
	
	Let us return to the proof of the Lemma \ref{lms}.
	Above  (Lemmas \ref{lj0}, \ref{lj}) it was shown that the Plucker and the Jacobi relations  belong to the ideal  
	$I_{g_n}$  of all relations between minors for the group  $G$.  Let us prove that these relations generate the ideal   $I_{g_n}$.
	
Indeed consider that mapping 	 $\varphi: X\mapsto \{a_X\}$, that maps a matrix to collection of it's minors.  Then the ideal  $I_{gl_m}$ is an ideal in the ring of polynomials in independent variables  $A_X$,  such that it's  null-space  is a closure of the image  $\varphi(GL_m)$.    When one passes from the image to it's closure only the origin of the point with zero coordinates is added.  Indeed let us do in each collection of coordinates  $\{a_X\}_{|X|=k}$  a projectivization  (separately in each collection).  The image   $\varphi(GL_m)$ after these projectivizations coincides with  a variety defined by Plucker relations   \cite{rel}.  Due to homogeneity of the Plucker relations and   homogeneity   for each  $k$ of the mapping  $pr_{\{a_X\}_{|X|=k}}\circ \varphi$ (where $pr_{\{a_X\}_{|X|=k}}$  is a projection to the written coordinates),  one obtains the needed statement.

	Thus if a collection of numbers  $\{a_X\}$ is non-zero and satisfy the Plucker relations than these numbers are minors of an invertible  $m\times m$ matrix. 
	
	Take  $m=2n$ or   $2n+1$  and consider an embedding  $G\subset GL_m$. The  ideal $I_{g_n}$ is an ideal whose null-space is a closure of  the image  $\varphi(G)$.
	
	Consider the ideal generated  by the Plucker and the Jacobi relations. Let   $\{a_X\}$  be a non-zero element from it null-space.
	 The Plucker relations provide that  $\{a_X\}$ can be considered as a collection of  minors of some matrix   $O$.
	
	Corollary  \ref{sl1} shows that  the Jacobi relations provide that   $\pm \Omega O^t\Omega=O^{-1}T$ for an upper-unitriangular matrix $T$, where $\pm=-$ for the series  $C$ and  $+$ for the series  $B$, $D$.
		This equality is equivalent to the following one:  $\pm O\Omega O^t\Omega=T $, from here one gets   that $T\Omega=O\Omega O^t$ is a skew-symmetric matrix for the series   $C$ and a symmetric matrix for the series $B$, $D$. Thus there exists a low-unitrangular matrix  $X$, such that $X^{}T \Omega X^{t}=\Omega$.
		
		For the matrix $XO$ one has  $\pm \Omega(XO)^t\Omega=(XO)^{-1}$, that is  $XO$  belongs to the considered group    $G$.  Since  $X$  is a low-unitrangular matrix then the determinants of type   \eqref{dete} for matrices  $XO$ and  $O$ coincide.
		
		Thus the relations from the ideal  $I_{\mathfrak{gl}_m}$ and the Jacobi relations provide that   $a_X$  are minors of a matrix   $XO\in G$.  As we have noted this implies that  the ideal generated by  $I_{gl_n}$ and the Jacobi relations is the ideal of all relations between the determinants   $a_X$ for the group  $G$.

		The case of the series  $D$  and determinants  $\bar{a}_X$  is reduced to the considered case by reordering of rows.   
		
	\end{proof}
	
	\subsubsection{ Some useful relations}

	Let us derive some relations that are corollaries of the Plucker and the Jacobi  relations that we shall use below. 
	
	\begin{propos}
		\label{dopson}
		In the case of the series $B$ one has
		
		\begin{equation}
			\label{s22}
			a_{\pm n,...,\widehat{\pm i}...,\pm 1,0}^{2}=-2\cdot a_{\pm n,...,\widehat{\pm i}...,\pm 1,-i}a_{\pm n,...,\widehat{\pm i}...,\pm 1,i}
		\end{equation}
		The choice of the sign on both sides is  the same.
		
	\end{propos}
	\begin{proof}
		
		Using a group automorphism that acts as a permutation of coordinates  $1,...,n$, one concludes that one can consider just the case   $i=1$
		
		One has the Plucker relations
		
		\begin{equation}
			\label{s211}
			a_{\pm n,...,\pm 2,-1}a_{\pm n,...,\pm 2,1,0}+a_{\pm n,...,\pm 2,1}a_{\pm n,...,\pm 2,-1,0}+a_{\pm n,...,\pm 2,,0}a_{\pm n,...,\pm 2,,-1,1}=0
		\end{equation}
		
		Take the Jacobi relations 
		
		\begin{align*}
			&a_{\pm n,...,\pm 2,1,0}= a_{\pm n,...,\pm 2,-1},&  a_{\pm n,...,\pm 2,-1,0}= a_{\pm n,...,\pm 2,-1},&  a_{\pm n,...,\pm 2,-1,1}=- a_{\pm n,...,\pm 2,0}.
		\end{align*}
		
		Applying them in  \eqref{s211}  one obtains \eqref{s22} for  $i=1$.

 The case of arbitrary $i$ is considered analogously.  the signs that appear because of permutations of indices are reducing.

\end{proof}

\begin{propos}
In the case of the series  $B$, $D$ one has

\begin{equation}
	\label{s333}	
	a_{\pm n,...,\widehat{\pm i},...,\widehat{\pm j},...,\pm 1,-i,i}=\sqrt{a_{\pm n,...,\widehat{\pm i},...,\widehat{\pm j},...,\pm 1,-i,-j}a_{\pm n,...,\widehat{\pm i},...,\widehat{\pm j},...,\pm 1,i,j}} ,
\end{equation}

in the case of the series  $D$  there exists an analogous relations for determinants   $\bar{a}_X$.
\end{propos}

\begin{proof}
As before let us consider first the case $i=1$, $j=2$.
To make notations less cumbersome let us omit the first indices.  One has the following Plucker relation 

$$
a_{...,-2,2}a_{...,-2,-1,1}+  a_{...,-2,1}a_{...,-2,2,-1}+ a_{...,-2,-1}a_{...,-2,1,2} =0,
$$

thus

\begin{align*}
	&  a_{...,-2,2}=-\frac{   a_{...,-2,1}a_{...,-2,2,-1}+ a_{...,-2,-1}a_{...,-2,1,2}  }{a_{...,-2,-1,1}}
\end{align*}

In the case of the series   $D$  one can embed  the Lie group into the  Lie group of the series   $B$  of  greater by one   dimension.  So we can assume that we deal with  the series   $B$. Now apply the  Jacobi relation to the determinants   $a_{...,-2,2,-1}, a_{...,-2,1,2}$ in the numerator. One obtains  

$$a_{...,-2,2}==\frac{   a_{...,-2,1}a_{...,0,-1}}{a_{...,-2,0}}+\frac{ a_{...,-2,-1}a_{...,0,1}  }{a_{...,-2,0}}$$

Now apply the relations  \eqref{s22}, one has

\begin{align*}
	&  a_{...,-2,2}=\frac{   a_{...,-2,1}\sqrt{a_{...,2,-1}a_{...,-2,-1}}}{\sqrt{a_{...,-2,1}a_{...,-2,-1}}}+\frac{ a_{...,-2,-1}\sqrt{a_{...,-2,1}a_{...,2,1}}  }{\sqrt{a_{...,-2,1}a_{...,-2,-1}}}=\\&=\sqrt{a_{...,-2,1}a_{...,2,-1}}+\sqrt{a_{...,-2,-1}a_{...,2,1}}
\end{align*}

Analogously the  exist the Plucker relations 

$$
a_{...,-1,1}a_{...,-1,-2,2}+a_{...,-1,2}a_{...,-1,1,-2}+a_{...,-1,-2}a_{...,-1,2,1}=0
$$

From here as above one gets 

\begin{align*}
	& a_{...,-1,1}=-\frac{a_{...,-1,2}a_{...,-1,1,-2}+a_{...,-1,-2}a_{...,-1,2,1}}{a_{...,-1,-2,2}}=
	\frac{a_{...,-1,2}a_{...,-2,0}}{a_{...,0,-1}}-\frac{a_{...,-1,-2}a_{...,0,2}}{a_{...,-1,0}}
\end{align*}

Now apply the relations  \eqref{s22}, one gets

\begin{align*}
	& a_{...,-1,1}=
	\frac{a_{...,-1,2}\sqrt{a_{...,-2,-1}a_{...,-2,1}}}{\sqrt{a_{...,-2,-1}  a_{...,2,-1}} }+\frac{a_{...,-1,-2}\sqrt{a_{...,-1,2}  a_{...,1,2} } }{\sqrt{ a_{...,-1,2} a_{...,-1,-2}    }}=\\&=
	-\sqrt{a_{...,2,-1}a_{...-2,1}}+\sqrt{a_{...,-1,-2}a_{...,1,2}}
\end{align*}

So one obtains

$$
a_{...,-2,2}+a_{...,-1,1}=2\sqrt{a_{...,-1,-2}a_{...,1,2}}
$$

From the Jacobi relations one gets that   $a_{...,-1,1}=a_{....,-2,2}$, so one gets

\begin{equation}
	a_{...,-1,1}=\sqrt{a_{...,-1,-2}a_{...,1,2}}=\sqrt{a_{...,-2,-1}a_{...,2,1}}
\end{equation}
The equality   \eqref{s333} for  $i=1$, $j=2$ is proved. The case of arbitrary  $i,j$  is considered analogously. The signs that appear because of permutations of indices  are reduced.

\end{proof}

	\subsection{ The conditions defining an irreducible representation in the Zhelobenko's model }
	\label{rzd3}
	
	\subsubsection{The general theorem}
	
	The indicators system is the system of differential equations of the following form:
	\begin{align}
		\begin{split}
			\label{indsys}
			&L_{-n,-n+1}^{r_{-n}+1}f=0,...,L_{-2,-1}^{r_{-2}+1}f=0,\,\,\,,L_{-1,1}^{r_{-1}+1}f=0\text{in the case of the series $A$, $C$},\\
			&L_{-n,-n+1}^{r_{-n}+1}f=0,...,L_{-2,-1}^{r_{-2}+1}f=0,\,\,\,,L_{-1,0}^{r_{-1}+1}f=0\text{ in the case of the series  $B$,}\\
			&L_{-n,-n+1}^{r_{-n}+1}f=0,...,L_{-2,-1}^{r_{-2}+1}f=0,\,\,\,L_{-2,1}^{r_{-1}+1}f=0\text{
				in the case of the series  $D$}.
		\end{split}
	\end{align}
	
	Here  $L_{-i,-j}$ is an operator acting onto a function  $f(a)$
	as a left infinitesimal shift by 
	$F_{-i,-j}$ for the series $B$, $C$, $D$ and by $E_{-i,-j}$ for the series $A$.

	The exponent  $r_{-i}$  is defined as follows:

	\begin{align}
		\begin{split}
			\label{rb}
			&r_{-n}=m_{-n}-m_{-n+1},...,r_{-2}=m_{-2}-m_{-1},\,\, r_{-1}=m_{-1}\text{in the case of the series  $A$, $C$,}\\
			&r_{-n}=m_{-n}-m_{-n+1},...,r_{-2}=m_{-2}-m_{-1},\,\, r_{-1}=2m_{-1}\text{ in the case of the series  $B$,}\\
			&r_{-n}=m_{-n}-m_{-n+1},...,r_{-2}=m_{-2}-|m_{-1}|,\,\,
			r_{-1}=m_{-2}+|m_{-1}|\text{ in the case of the series $D$}.
		\end{split}
	\end{align}

	In \cite{zh}  the following was proved.
	\begin{theorem}\label{tzh0}
		For the series  $A$ one has the following statement. Let the group   $GL_{n+1}$ act in the space with coordinates indexed by   $-n,...,-1,1$\footnote{such an indexation is used for a formulation of analogous results for other series}. Then in  the space of functions on the group   $GL_{n+1}$  an irreducible representation with the highest weight  $[m_{-n},...,m_{-1},0]$ and the highest vector  \eqref{stv0}
		is defined by the following conditions.

		\begin{enumerate}
			\item $L_{-}f=0$, where  $L_{-}$ is a left infinitesimal shift by an arbitrary element of $GL_{n+1}$ corresponding to a negative root.
			\item $L_{-i,-i}f=m_{-i}f$, where $L_{-i,-i}$, $i=-1,1,...,n$ is a left infinitesimal shift by an element of the group 
			$GL_{n+1}$, corresponding to a Cartan element    $E_{-i,-i}$.
			\item $f$  satisfies the indicator system.
		\end{enumerate}
		
	\end{theorem}
	
	This  Theorem describes the space of the model $Zh$ in the case of the series $A$.
	Let us prove an analogous statement for the series  $B$, $C$, $D$.

	\begin{theorem}
		\label{tzh}
		In the Zhelobenko's realization for the series  $B$, $C$, $D$ an irreducible representation with the highest vector  \eqref{stv0}
		is defined by the conditions $1,3$  and the condition $2$ where $E_{-i,-i}$ is changed to	$F_{-i,-i}$ in the definition of   $L_{-i,-i}$.

	\end{theorem}
	
	The  proceeding part of  Section   \ref{rzd3} is devoted to the proof of this Theorem.
	
	\begin{proof}
		
		The scheme of the proof is the following. First of all one gives formulas for the action of an operator of a left infinitesimal shift onto determinant. Using them one proves that the highest vector satisfies the conditions 1-3. The main difficult is to prove that the highest vector satisfies the indicator system. Then one easily proves that an arbitrary vector satisfies the conditions  1-3. Finally one proves that among  the functions that satisfy the conditions   1-3 there is nothing but functions that vectors of the representation with the highest vector 	\eqref{stv0}.

		Let us proceed to realization of this plan.	
		Consider the determinant
		$a_{i_1,...,i_k}$
		and let us find an action of the left infinitesimal shift onto it. Let us introduce a temporary notation for determinants
		
		$$
		a^{-n,...,-n+k-1}_{i_1,...,i_k},
		$$
		where we write also the indices of rows (on the up). Then the operator   $L_{-i,-j}$ of the left infinitesimal shift acts onto row indices $-n,...,-n+k-1$ by the following ruler. For the series  $A$  the operator of the left infinitesimal shift $E_{-i,-j}$   acts
		by the ruler
		
		\begin{equation}
			\label{lij}
			L_{-i,-j}a^{-n,...,-n+k-1}_{i_1,...,i_k}=a^{\{-n,...,-n+k-1\}\mid_{-i\mapsto -j}}_{i_1,...,i_k}.
		\end{equation}
		
		The action  of the operator of the left infinitesimal shift by   $F_{-i,-j}$ for the series $B$,
		$C$, $D$  can be expressed through these operators.
		
		Analogously one proves that 
		
		\begin{equation}
			\label{fii}
			L_{-i,-i}a^{-n,...,-n+k-1}_{i_1,...,i_k}=\begin{cases}
				a^{-n,...,-n+k-1}_{i_1,...,i_k}\text{ if  }
				-i\in\{-n,...,-n+k-1\},\,\,\, \\ \text{ $0$ otherwise }\end{cases}
		\end{equation}
		
		Onto a product of determinants $L_{-i,-i}, L_{-i,-j}$ act by the Leibniz ruler.
		
		Let us prove that the conditions 1-3 for a representation with the highest weight  \eqref{stv0}  are satisfied.
		
		Let us prove that the conditions  1-3  are satisfied for the highest vector  \eqref{stv0}.  The fact that the conditions  1 and 2 are satisfied is proved directly. One has to prove that the condition 3 is satisfied.
		
		From the formula  \eqref{lij} it follows that the operators  $L_{-i,-i+1}$ for $i=-n,...,-3$  being acting onto  \eqref{stv0}  actually are acting onto  $a_{-n,...,-i}^{m_{-i}-m_{-i+1}}$ only.  The operator $L_{-2,-1}$ for the series $A$, $B$, $C$ and  $D$  in the case  $m_{-1}\geq 0$,  the operator $L_{-2,1}$
		for the series  $D$ in the case  $m_{-1}<0$
		being acting onto  \eqref{stv} is acting onto $a_{-n,...,-2}^{m_{-2}-m_{-1}}$ only.
		The operator  $L_{-1,1}$ for the series  $A$, $C$,  the operator $L_{-1,0}$  for the series  $B$ are acting onto determinants of order  $n$ only.
		
		In the case of the series   $D$  one also has the following operators. In the case $m_{-1}\geq 0$ one has the operator $L_{-2,1}$,
		and in the case  $m_{-1}< 0$ -  one has an operator  $L_{-2,-1}$. Both they act onto determinants of orders   $n-1$ and $n$.
		
		Since the exponents  $r_{-i}=m_{-i}-m_{-i+1}$ of the determinants  $a_{-n,...,-n+i-1}$ are integer then by the formula  \eqref{lij}  the equation    $L_{-i,-i+1}^{r_{-i}+1}v_0=0$, $i=-n,...,-2$  is satisfied.  By the same reason the equation  $L_{-1,1}^{r_{-1}+1}v_0=0$  is satisfied for the series  $A$, $C$.

		Let us check that equality  $L_{-1,0}^{m_{-1}+1}v_0=0$ holds for the series  $B$.
		To prove it one needs the following formulas. In the case of the series   $B$ one has

		\begin{align*}
			&L_{-1,0}a_{-n,...,-2,-1}^{-n,...,-2,-1}=a_{-n,...,-2,-1}^{-n,...,-2,0},\,\,\,L_{-1,0}^2a_{-n,...,-2,-1}^{-n,...,-2,-1}=-a_{-n,...,-2,-1}^{-n,...,-2,1},\\&L_{-1,0}^3a_{-n,...,-2,-1}^{-n,...,-2,-1}=0.
		\end{align*}
		
		Using these formula one derives that the equation $L_{-1,0}^{2m_{-1}+1}v_0=0$ holds in the case of the series  $B$ when the highest weight in integer.
		
		Also one has

		\begin{align*}
			&L_{-1,0}(a_{-n,...,-2,-1}^{-n,...,-2,-1})^{1/2}=\frac{1}{2}a_{-n,...,-2,-1}^{-n,...,-2,0}(a_{-n,...,-2,-1}^{-n,...,-2,-1})^{-1/2},\\
			&L_{-1,0}^2a_{-n,...,-2,-1}^{-n,...,-2,-1}=-\frac{1}{2}a_{-n,...,-2,-1}^{-n,...,-2,1}(a_{-n,...,-2,-1}^{-n,...,-2,-1})^{-1/2}-\\&-
			\frac{1}{4}(a_{-n,...,-2,-1}^{-n,...,-2,0})^2(a_{-n,...,-2,-1}^{-n,...,-2,-1})^{-3/2}=0
		\end{align*}
		
		To obtain the last equality one uses a relation  $$a_{-n,...,-2,-1}^{-n,...,-2,1}a_{-n,...,-2,-1}^{-n,...,-2,-1}=-\frac{1}{2}(a_{-n,...,-2,-1}^{-n,...,-2,0})^2.$$

		From these equalities  it follows that $L_{-1,0}^{2m_{-1}+1}v_0=0$ in the case of the series  $B$  when the highest weight is half-integer.

		Now consider the case of the series  $D$ when $m_{-1}\geq 0$. Let us check the the equality $L_{-2,1}^{m_{-2}+m_{-1}+1}v_0=0$ holds. One has 
		\begin{align*}
			&L_{-2,1}a_{-n,...,-2}^{-n,...,-2}=a_{-n,...,-2}^{-n,...,1},\,\,\, L_{-2,1}a_{-n,...,-2}^{-n,...,1}=0,\\
			&L_{-2,1}a_{-n,...,-2,-1}^{-n,...,-2,-1}=2a_{-n,...,-2,-1}^{-n,...,1,-1},\,\,\, L_{-2,1}a_{-n,...,-2,-1}^{-n,...,1,-1}=a_{-n,...,-2,-1}^{-n,...,1,2},\,\,\, L_{-2,1}a_{-n,...,-2,-1}^{-n,...,1,2}=0.
		\end{align*}
		
		Thus  $v_0$ is mapped to zero under the action of   $L_{-2,1}$ to the power which equals to a sum of   $1$  and the size of the determinant 	$a_{-n,...,-2}$ and twice of the size of the determinant    $a_{-n,...,-2,-1}$. That is under the action of   $L_{-2,1}$  to the power 
		$1+(m_{-2}-m_{-1})+2m_{-1}=1+m_{-2}+m_{-1}$.
		
		In the case of the series  $D$ when $m_{-1}< 0$.  Let us check the the equality $L_{-2,-1}^{m_{-2}+m_{-1}+1}v_0=0$ holds.   One has
		
		\begin{align*}
			&L_{-2,-1}a_{-n,...,-2}^{-n,...,-2}=a_{-n,...,-2}^{-n,...,-1},\,\,\, L_{-2,-1}a_{-n,...,-2}^{-n,...,-1}=0,\\
			&L_{-2,-1}a_{-n,...,-2,1}^{-n,...,-2,1}=2a_{-n,...,-2,-1}^{-n,...,-1,1},\,\,\, L_{-2,-1}a_{-n,...,-2,-1}^{-n,...,-1,1}=-a_{-n,...,-2,-1}^{-n,...,-1,2},\,\,\, L_{-2,1}a_{-n,...,-2,-1}^{-n,...,-1,2}=0.
		\end{align*}
		
		Thus $v_0$ is mapped to zero under the action of    $L_{-2,-1}$ to the power
		$1+(m_{-2}-m_{-1})+2m_{-1}=1+m_{-2}+m_{-1}$.

		So  the highest vector satisfies the indicator system.
		
		The fact that the conditions  1-3 do hold for an arbitrary vector follows from the fact that left and right shifts commute and an arbitrary vector can be written as a linear combination of right shifts of the highest vector.

		Now one has to prove that among functions satisfying the conditions 1-3,  there is nothing  but functions that form a representation with the highest vector  \eqref{stv0}.

		Let us be given a function on   $G$, that satisfies the conditions 1-3.  Then it's restriction to the subgroup $Z\subset G$ of upper-unitriangular matrices satisfies the indicator system. According to  \cite{zh} this restriction belongs to a realization of the representation in the space of functions on the group   $Z$ (in \cite{zh} this realization in the space of functions on   $Z $ is given for all series).  Using the results of   \cite{zh}, one gets that the initial function on   $G$ belongs to the representation with the highest weight  \eqref{stv0}.

	\end{proof}

	\subsubsection{Solution of the indicator system and the equations $L_{-i,-i}f=m_{-i}f$}
	\label{islm}

	In the proof of the Theorem   \ref{tzh}  the formulas for the action of $L_{-i,-i}$ were obtained. One gets using them the following  Lemma.
	
	\begin{lemma}
		\label{lfmf}
		Solution of the system $L_{-i,-i}f=m_{-i}f$ that are functions of determinants are described as follows.  If one represents this function as  a sum of monomials in determinants  then in each monomial the sum of exponents of determinant of size  $n-i+1$ equals
		$r_{-i}$ for $i=n,...,2$.  Also the sum of exponents of determinants of size $n$ equals $|m_{-1}|$.

	\end{lemma}

	Now has to find solutions of the indicator system among these functions.
	
	\begin{lemma}
		\label{l1} If the highest weight is integer non-negative then the solutions of the indicator system are polynomials in determinants that satisfy the conditions of the Lemma   \ref{lfmf}.
	\end{lemma}
	
	\begin{proof}
		Let highest vector
		\eqref{stv0} be a polynomial in determinants. Since under the action of the algebra  of polynomials in determinant is invariant then an arbitrary vector of the representation is also presented as a polynomial in determinants. Using this argument one immediately derives the statement of the Lemma from   \eqref{lij}.
	\end{proof}

	\begin{lemma}
		\label{l2} If the highest vector in half-integer then among the functions that satisfy the conditions of the Lemma 
		\ref{lfmf}, the solutions of the indicator system are the functions of type   
		
		\begin{align}
			\begin{split}
				\label{f12}
				&f=\sum_{\alpha}\sqrt{a_{\pm n,...,\pm 2,\pm 1}}\cdot f_{\alpha}\text{ in the case of the series $B$  and of $D$ when $m_{-1}\geq 0$},\\
				&f=\sum_{\alpha}\sqrt{\bar{a}_{\pm n,...,\pm 2,\pm 1}}\cdot f_{\alpha}\text{
					in the case of the series  $D$ and $m_{-1}< 0$}.
			\end{split}
		\end{align}
		here $\alpha$ denotes a choice of  $+$ or $-$  for all the indices  $\pm n,...,\pm 1 $.   In the case of the series  $B$ no conditions on the choice of the signs are imposed. I the case of the series   $D$  when $m_{-1}\geq 0$  the sign  $+$ is chosen in such  way that the parity of the quantity of  places with $-$ is the same as the parity of $n$,  and in the case  $m_{-1}< 0$ the parity  of the quantity of  places with $-$ must be the same as the parity of   $n-1$.
		
		All $f_{\alpha}$  are polynomials in determinants and the sum of exponents of determinants of size  $n-i+1$ equals
		$r_{-i}$ for $i=n,...,2$. In the case of the series  $D$ and $m_{-1}<0$  one takes as determinants of order   $n$ the determinants of type   $\bar{a}_X$.	
		In all cases the sum of exponents of  determinants order $n$ equals  $|m_{-1}|-\frac{1}{2}$.
	\end{lemma}
	
	\begin{proof}
		Consider the case $B$.
		
		First of all let us show that functions representing vectors of the representation with the highest weight \eqref{stv0} have the form	\eqref{f12}.

		Consider first a representation with the highest vector  $\sqrt{a_{-n,...,-2,-1}}$,  that is the spinor representation

		\begin{propos}
			\label{spb} For the algebra $\mathfrak{o}_{2n+1}$
			a representation with the highest vector  $\sqrt{a_{-n,...,-2,-1}}$ is a span of functions of type  
			
			\begin{equation}
				\label{spn}
				\sqrt{a_{\pm n,...,\pm 2,\pm 1}}
			\end{equation}

		\end{propos}
		
	\begin{proof}
Let us prove that when one applies 	  $F_{p,q}$ to  $
	\sqrt{a_{\pm n,...,\pm 2,\pm 1}}
	$ one obtains a linear combination of function of the same type. Then the considered linear span is a representation that contains the spinor representation since   $\sqrt{a_{-n,...,-2,-1}}$ is the higest vector of the spinor representation.    The coincidence of the considered linear span and the spinor representation follows from the fact that both these linear spaces have dimension   $2^n$.

The proof the statement formulated in the paragraph above is the same for choices of signs.  To avoid cumbersome notations we assume that 	 everywhere the sign  $-$ is chosen.  Then a non-zero result is obtained under the action onto   $\sqrt{a_{-n,...,-2,-1}}$  only for generators   $F_{0,-i}$, $F_{j,-i}$, $i,j>0$.   In these cases one gets 
	
	\begin{align}
		\begin{split}
			\label{ffa}
			&F_{0,-i}:\,\,\sqrt{a_{-n,...,-2,-1}} \mapsto\frac{(-1)^i a_{-n,...,\widehat{-i},...,-1,0}}{2\sqrt{a_{-n,...,-1}}},\\
			&F_{j,-i}:\,\,\sqrt{a_{-n,...,-2,-1}} \mapsto\frac{(-1)^i a_{-n,...,\widehat{-i},...,-1,j}-(-1)^j   a_{-n,...,\widehat{-j},...,-1,i} }{2\sqrt{a_{-n,...,-1}}}.
		\end{split}
	\end{align}

	Consider the first equality in  \eqref{ffa}.
	One has a relation \eqref{s22}. Taking a square root from it one gets 
	
	\begin{equation}
		\label{fl1}
		F_{0,-i}\sqrt{a_{-n,...,-2,-1}}=\frac{\sqrt{-2}}{2} \cdot \sqrt{a_{-n,...,\widehat{-i},...,-1,i}}
	\end{equation}
	
	
	Now take the second equality in   \eqref{ffa}.  Let us prove that 
	$$
	\frac{(-1)^i a_{-n,...,\widehat{-i},...,-1,j}-(-1)^j   a_{-n,...,\widehat{-j},...,-1,i}}{2\sqrt{a_{-n,...,-1}}}
	$$
is a linear combination of functions 	  \eqref{spn}.  Without loss of generality one can assume that    $i=1$, $j=2$.  That is we consider the fraction 
	
	\begin{equation}
		\label{drb}
		\frac{  a_{-n,...,-3,-2,2}- a_{-n,...,-3,1,-1}}{2\sqrt{a_{-n,...,-1}}}
	\end{equation}

	Using  \eqref{s333} (for   $i=1$,  $j=2$ and expressing  $a_{-n,...,-3,-1,1}$, and also for   $i=2$,  $j=1$  and expressing $a_{-n,...,-3,-2,2}$) one gets that  \eqref{drb} equals $\sqrt{a_{...,2,1}}$. In other words when one substitutes into  \eqref{ffa}, for  $i=1,j=2$ one gets
	
	\begin{equation}
		\label{chtn}
		F_{2,-1}\sqrt{a_{-n,...,-1}}= \sqrt{a_{-n,...,2,1}}.
	\end{equation}
	
\end{proof}

		Now consider the case of an arbitrary highest weight. The  highest vector
		can be represented as
		
		$$
		v_0=v'_0(a_{-n,...,-2,-1})^{\frac{1}{2}},
		$$
		where  $v'_0$  is a polynomial in determinants. An arbitrary vector $f$
		can be written as a linear combination of vectors that are obtained by the action onto the highest vector of the operators $\prod_{i<j} F_{i,j}^{p_{i,j}}$.  As a result one obtains  a vector of type 
		\eqref{f12}.

		Thus the Zhelobenko's model is contained in the space of functions of type  \eqref{f12}.   Now one has to prove that all function of this type belong to the  Zhelobenko's model. For this it is enough to check that each function of the form \eqref{f12},
		that satisfies the conditions of Lemma \ref{lfmf}, is a solution of the indicator system. 
		
		The operators $L_{-i,-i+1}$, $i=n,...,2$  are acting onto the determinant of size  $n-i+1$ only. Such determinants occur in   $f_{\alpha}$ only, thus their  exponents are non-negative integers. The sum of the exponents of determinants of size  $i$ equals  $r_{-n+i-1}$. Thus the conditions of Lemma \ref{lfmf} are satisfied.  Thus the conditions 
		$L_{-i,-i+1}^{r_{-i}+1}f=0 $ for $i=n,...,2$ hold.
		
		Now consider the equation
		$L_{-1,0}^{2m_{-1}+1}f=L_{-1,0}^{2[m_{-1}]+1+1}f=0$,  where
		$[m_{-1}]$ is the integer part.  The operator $L_{-1,0}^{2[m_{-1}]+1+1}$ acts by the Leibniz ruler onto each summand in 
		\eqref{f12} as follows.
		
		Either  $L_{-1,0}^{2[m_{-1}]+1+1}$ acts onto the second factor  $f_{\alpha}$. In this case one obtains    $0$, since in the polynomials  $f_{\alpha}$ the sum of exponents of determinants of size $n$ equals  $[m_{-1}]$,  but such polynomials are annihilated by$L_{-1,0}^{2[m_{-1}]+1}$.
		
		Either   $L_{-1,0}^{2[m_{-1}]+1}$ acts onto the second factor   $f_{\alpha}$,  and 
		$L_{-1,0}$ acts onto the first factor.  In this case one obtains  $0$  by the same reason.
		
		Either $L_{-1,0}^{2[m_{-1}]+2-k}$ acts onto the second factor   and  $L_{-1,0}^{k}$  acts onto the first factor, here $k\geq 2$.  Since the first factor is a vector of the representation with the highest weight  $[\frac{1}{2},...,\frac{1}{2}]$, then under the action of  $L_{-1,0}^{2}$
		it vanishes.
		
		Thus the vectors of type \eqref{f12}  are vanishing under the action of 
		$L_{-1,0}^{2m_{-1}+1}$.  In the case of the series $B$ the Lemma is proved.
		
		Now consider the case of the series  $D$.   The scheme of the proof is the same as in the case of the series    $B$.  First of all one considers the spinor representations with the highest weights  $[\frac{1}{2},...,\frac{1}{2},\frac{1}{2}]$,   $[\frac{1}{2},...,\frac{1}{2},-\frac{1}{2}]$.

		\begin{propos}
			For the algebra  $\mathfrak{o}_{2n}$ the representations with the highest vectors    $\sqrt{a_{-n,...,-2,-1}}$, $\sqrt{\bar{a}_{-n,...,-2,1}}$  coincide with the spans of functions 
			
			\begin{align*}
				&  < \sqrt{a_{\pm n,...,\pm 2,\pm1}}> ,\text{ the parity of the number of $-$ equals to the parity of $n$},\\
				&  < \sqrt{\bar{a}_{\pm n,...,\pm 2,\pm1}}>, \text{ the parity of the number of $-$ equals to the parity of $n-1$}
			\end{align*}

		\end{propos}
		
		\begin{proof}
			
			The scheme of the proof is the same as for the Proposition \ref{spb}.
			One has an embedding of the Lie algebras $\mathfrak{o}_{2n}\subset\mathfrak{o}_{2n+1}$, induced by an embedding of root systems.  Thus the determinants  $a_{\pm n,...,\pm 1}$ and also  $\bar{a}_{\pm n,...,\pm 1}$ can be considered as functions on the group $SO_{2n+1}$. Thus we can use the formula  \eqref{chtn}.
			
			Using it one obtains that under the action of $\mathfrak{o}_{2n}$  onto a function $\sqrt{a_{\pm n,...,\pm 2,\pm1}}, \sqrt{\bar{a}_{\pm n,...,\pm 2,\pm1}} $ one obtains a linear combination of functions of the same type. At the same time, the parity of the number of minuses in the indexes is preserved. Thus, both linear spans under consideration are representations. The dimensions of each of these linear spans are equal  $2^{n-1}$.

			The irreducible representations with the highest vectors $\sqrt{a_{-n,...,-2,-1}}$, $\sqrt{\bar{a}_{-n,...,-2,1}}$ also have dimension $2^{n-1}$.
			
			From these two facts one obtains the statement of the Proposition.

		\end{proof}
		
		Further considerations in the case of the series  $D$  are analogous to considerations in the case of the series   $B$. The only change is a replacement of $L_{-1,0}$ onto  $L_{-2,1}$.
		
	\end{proof}

	Thus we have proved the following Theorem.
	
	\begin{theorem}
		\label{ost1}
		The space of an irreducible representation with the highest vector \eqref{stv0} in the Zhelobenko's model is described as follows.
		In the cases of the series  $A$, $B$, $C$   and the series $D$ with $m_{-1}\geq 0$    consider the space of functions of determinants   $a_{X}$, $|X|\leq n$.  And  in the case of the series $D$ with  $|m_{-1}|<0$  consider the space of functions of determinants  $a_{X}$, $|X|< n$  and $\bar{a}_X$, $|X|=n$.
		
		The functions must satisfy  the following condition. In the case of an integer of highest weight, they are polynomials in the determinants satisfying the conditions of the Lemma \ref{lfmf}.
		In the case of an integer of highest weight they are function of type \eqref{f12},  satisfying the conditions of the Lemma   \ref{l2}.
	\end{theorem}

	\section{ The GKZ system for the series  $A$. The A-GKZ system}
	\label{seriaa}
	In this s section  we introduce two important systems of partial differential equations. Solutions of one of the first  system (the GKZ system) are used to form a basis in the Zhelobenko model, and solutions of the second (the A-GKZ system) are used to build a new representation model.

	The results of this section concern only the $A$ series, essentially all of them were obtained in \cite{a4}, but we modify them insignificantly to adapt them for consideration of other series. It is shown how proofs of modified statements can be obtained from proofs of similar statements in \cite{a4}.

	\subsection{The Gelfand-Tsetlin lattice. The vectors $v_{\alpha}$. The number $\mathcal{K}$}
	\label{rgc}

	Consider the Lie algebra $\mathfrak{gl}_m$, which we identify with the algebra of matrices whose rows and columns are indexed by the numbers $-n,...,\hat{0},...,n$ (in the case $m=2n$) or $-n,...,0,....n$ (in the case $m=2n+1$). Le us associate with it a shifted lattice in the space $\mathbb{Z}^{N}$, whose coordinates are numbered by proper subsets $X\subset \{-n,...,n\}$. Here $N$ is the number of possible proper subsets of $X\subset\{-n,...,n\}$. Such a strange at first glance indexing is taken in order to obtain a lattice that will be used in further reasoning for the series $B$, $C$, $D$.

	The definition of the Gelfand-Tsetlin lattice for the $A$ series, which is given below, differs somewhat from that given in \cite{a4}. The essence of the difference is as follows. The idea of a homogeneous system defining the Gelfand-Tsetlin lattice is that an inhomogeneous versions of this system represent naive conditions on the vector of exponents of a monomial  in  determinants, which can appear in the decomposition of a function corresponding to the vector of the Gelfand-Tsetlin basis.
	Accordingly, this system should depend on the chain of subalgebras underlying the construction of the Gelfand-Tsetlin basis. In \cite{a4}, a chain is used in which the subalgebra $\mathfrak{gl}_{m-k}$ is formed by matrices having nonzero elements of the first $m-k$ rows and columns relative to the standard order on the set $\{-n,...,n\}$. In the present paper, we use a chain in which the subalgebra $\mathfrak{gl}_{m-k}$ is formed by matrices having nonzero elements in $m-k$ rows and columns with the smallest indices relative to the ordering

	\begin{equation}\label{porya}1\succ -1\succ 2\succ -2\succ...\succ 0 ,\end{equation} 
	
	This is done in order to use the results for the $A$ series in further consideration of other series.

	We also introduce the value

	\begin{equation}
		\label{snk}
		s(p,q):=\#\{t:t\leq p\,\, \& \,\,t\succ q\}     
	\end{equation}

	\begin{definition}
		Define the Gelfand-Tselin lattice $\mathcal{B}_{GC}^{\mathfrak{gl}_m}\subset \mathbb{Z}^N$ for the series   $A$ by the following homogeneous system of equations

		\begin{align}
			\begin{split}
				\label{sar}
				&\delta\in \mathcal{B}_{GC}^{\mathfrak{gl}_m} \Leftrightarrow  \forall p,q \in \{-n,...,n\},p\preceq q  \text{ one has  }:\\&\sum_{X\text{ contains }\geq (p+n+1-s(p,q))\text{ elements }\preceq  q} \delta_{X}=0.
			\end{split}
		\end{align}

		The resulting lattice is called the Gelfand-Tsetlin lattice for the series $A$.
	\end{definition}
	
	Let us construct the generators of this lattice. Take subsets  $Y_{i},i=-n,...,n$  of the following type
	
	\begin{align}
		\begin{split}
			\label{yi}
			&Y_{-n}=\{-n\},\,\,\, Y_{n}=\{-n,n\},\\
			&...\\
			& Y_{-i}=\{-n,...,-i\},\,\,\,\, Y_{i}=\{-n,...,-i,i\},\,\,\, i>0
		\end{split}
	\end{align}
	
	Note that  $Y_{-i}=Y_{-i-1}\cap\{-i\}$, $Y_{i}=Y_{-i}\cap\{i\}$.  Also note that 
	
	$$
	\{j: j\succ Y_{-i}\}=\{j:j\succeq i\},\,\,\, \{j: j\succ Y_{i}\}=\{j:j\succeq -i+1\},
	$$
	
	where the index is greater than the subset if  it is greater than each element of the subset.

	Define the vectors (in the case  $i\geq 0$)
	
	\begin{align}
		\begin{split}
			\label{va0}
			&v_{\alpha}=e_{Y_{-i-1},-i,X}-e_{Y_{-i-1},j,X}-e_{Y_{-i-1},-i,y,X}
			+e_{Y_{-i-1},j,y,X}\text{ for }-i\prec j\prec y \prec X\\
			&v_{\alpha}=e_{Y_{-i},i,X}-e_{Y_{-i},j,X}-e_{Y_{-i},i,y,X}
			+e_{Y_{-i},j,y,X}\text{ for }i\prec j\prec y \prec X
		\end{split}
	\end{align}
	
	Here  $e_{Z}$  is a unit base vector corresponding to a coordinate  $Z$.  One  checks directly  that these vectors satisfy the system   \eqref{sar}.

	To select vectors, we fix $\pm i$, $j$, $y$ and consider vectors \eqref{va0} for all possible $X$. Let's construct a graph in which these subsets  $X$ will be the vertices, and edges will be pairs of subsets of the form $X_1=X, X_2=\{x\}\cup X$. Let's choose such a collection of subsets that the corresponding subgraph is a tree, and it is  the maximum  collection with respect to the extension while preserving this property.	For fixed $\pm i$, $j$, $y$ and chosen $X$, we construct vectors \eqref{va}. 
	
	
	
	\begin{lemma}\label{lm}
		The chosen vectors  $v_{\alpha}$  form a base in  the lattice  $\mathcal{B}_{GC}^{\mathfrak{gl}_m}$.
	\end{lemma}
	
The proof of this Lemma is analogous to the proof of an analogous statement in  \cite{a4} and that is why it is omitted here.

	 Let us write  (for  $i\in \{-n,...,n\}$)
	
	\begin{equation}
		\label{va}
		v_{\alpha}=e_{i,Z}-e_{j,Z}-e_{i,y,Z}
		+e_{j,y,Z},\end{equation}
	
		where  $Z=Y_{i-1}\cup X$
		for $i\leq 0$ and  $Z=Y_{-i}\cup X$
		for $i >0$.

	Denote as  $\mathcal{K}$ the number of vectors in a base of $\mathcal{B}_{GC}^{\mathfrak{gl}_m}$.

	\subsection{The GKZ system for the series  $A$}
	
	Let's write  the GKZ system associated with the lattice $\mathcal{B}_{GC}^{\mathfrak{gl}_m}$. In the Section \ref{r2}, a GKZ system associated with an arbitrary lattice was defined. In this section we will write only equations of the second type, which are determined by the lattice $\mathcal{B}_{GC}^{\mathfrak{gl}_m}$. This will be the GKZ system for the $A$ series.
	
	Namely, let $A_X$ be variables that are skew-symmetric over the set of $X$, but they do not obey any relations.
	Then we  call the GKZ system for the $A$ series a system {\it generated }  by differential equations constructed from vectors \eqref{va}
	
	\begin{equation}
		\label{gkza1}
		(\frac{\partial^2}{\partial A_{i,Z}\partial A_{j,y,Z} }-\frac{\partial^2}{\partial A_{j,Z}\partial A_{i,y,Z} })\mathcal{F}
	\end{equation}

Below when we refer  to the system  \eqref{gkza1}  we mean the system generated by these equations. This  agreement  is used for all GKZ and A-GKZ system  introduced in the present paper. 
	
	\begin{lemma}[\cite{GG}]
		\label{lgkz}
		In the space of polynomial solution of the system   \eqref{gkza1} there exists a base $\mathcal{F}_{\gamma}(A;\mathcal{B}_{GC}^{\mathfrak{gl}_m})$, consisting of   $\Gamma$-series.  We take that shift vectors $\gamma$ that represent different elements of the factor-space  $\mathbb{C}^N /\mathcal{B}_{GC}^{\mathfrak{gl}_m}$ and such that one of the representatives of this class has only non-negative coordinates  (see Section   \ref{dgca}).
	\end{lemma}
	
	In the formulation, the term <<polynomial>> can be replaced by <<decomposing into power series>>, but it should be understood that the system has other solutions decomposing into power-logarithmic series.
	
	\subsection{The A-GKZ system for the series  $A$}
	\label{agkzs}
	
	Let us call the   A-GKZ system for the series   $A$   the system  {\it generated }  by following equations:
	
	\begin{equation}
		\label{gkza}
		(\frac{\partial^2}{\partial A_{i,Z }\partial A_{j,y,Z} }-\frac{\partial^2}{\partial A_{j,Z}\partial A_{i,y,Z} }+\frac{\partial^2}{\partial A_{y,Z}\partial A_{i,j,Z} })F
	\end{equation}
	\begin{remark}
		
		An important observation is that with an equation of the A-GKZ system (i.e. with the vector \eqref{va}))  a Plucker relation for determinants is associated 
		\begin{equation}
			\label{spgl}
			a_{i,Z} a_{y,Z}	- a_{j,Z} a_{i,y,Z}	+a_{y,Z} a_{i,j,Z}	=0
		\end{equation}
	\end{remark}

	With the vector   $v_{\alpha}$ of type \eqref{va}  one relates a vector

	\begin{equation}
		\label{ra}
		r_{\alpha}=e_{y,Z}-e_{j,Z}-e_{i,y,Z}+e_{i,j,Z}
	\end{equation}
	
	Now let us give formulas for some solutions of the A-GKZ system. Let $t,s\in\mathbb{Z}^\mathcal{K}_{\geq 0}$, where $\mathcal{K}$ be the number of independent vectors $v_{\alpha}$. Introduce the functions

	\begin{equation}
		\label{fs}
		\mathcal{F}_{\gamma}^s(A;\mathcal{B}_{GC}^{\mathfrak{gl}_m})=\sum_{t\in\mathbb{Z}^{\mathcal{K}}}\frac{(t+1)...(t+s)A^{\gamma+tv}}{s!(\gamma+tv)!},
	\end{equation}
	in the formula a multi-index notation is used  $tv:=t_1v_1+...+t_{\mathcal{K}}v_{\mathcal{K}}$, $sr:=s_1r_1+...+s_{\mathcal{K}}r_{\mathcal{K}}$  and also a multi-index notation for factorials. Put

	\begin{equation}
		\label{ffpr}
		F_{\gamma}(A;\mathcal{B}_{GC}^{\mathfrak{gl}_m})=\sum_{s\in \mathbb{Z}^{\mathcal{K}}}(-1)^s\mathcal{F}_{\gamma-sr}^s(A;\mathcal{B}_{GC}^{\mathfrak{gl}_m}).
	\end{equation}
	
	Note that this series does depend on the choice of the vector  $\gamma$, not only on the class $\gamma\,\,mod \mathcal{B}_{GC}^{\mathfrak{gl}_m}$.
	As in   \cite{a4}  (see also  \cite{a5}) one proves
	
	\begin{propos}
		\label{sagkz}
		The functions  $F_{\gamma}(A;\mathcal{B}_{GC}^{\mathfrak{gl}_m})$, for vectors  $\gamma$ chosen in the Lemma \ref{lgkz}, form a base in the space of polynomial solutions of the A-GKZ system.
	\end{propos}

	\begin{remark}\label{rem1}

		In \cite{a4} the following statement is proved. It's  proof  is again without any change is valid in the considered in the present paper situation.  Let $a_{i_1,...,i_k}$ be a function of the form \eqref{dete} for the group $GL_m$. Let $I_{\mathfrak{gl}_m}$ be the ideal of relations between these minors, considered as an ideal in the polynomial ring $\mathbb{C}[A]$.

		To any relations their corresponds a differential operator by the ruler   \begin{equation}\label{aa}A_X\mapsto \frac{\partial}{\partial A_X}.\end{equation} Thus one gets an ideal   $\bar{I}_{\mathfrak{gl}_m}$
		in the ring of differential operators with constant coefficients.	
		t can be noted that the equations of the A-GKZ system are obtained by the \eqref{aa} rule from {\it of some} relations of the Plucker type.
		
		In  \cite{a4} it is verified that, nevertheless, the solution space of the ideal $\bar{I}_{\mathfrak{gl}_m}$ coincides with the solution space of the A-GKZ system.

	\end{remark}
	
	\section{ Models of representations for  $\mathfrak{gl}_m$}
	
	The results of this section also concern only the $A$ series, they were obtained in \cite{a4}.

	As before, let's take the Lie algebra $\mathfrak{gl}_m$, where $m=2n$ or $2n+1$, acting in a space with coordinates $-n,...,n$, where in the first case $0$ is thrown out, and in the second - not. Next, $\{-n,...,n\}$ denotes such an index set.

	\subsection{Gelfand-Tsetlin diagrams and shift vectors}
	\label{dgca}

	With a Gelfand-Tsetlin diagram  $(m_{p,q})$
	for the algebra  $\mathfrak{gl}_m$ one associates a shifted lattice $\Pi=\gamma+\mathcal{B}_{GC}^{\mathfrak{gl}_m}$, defined by an inhomogeneous system of equations

	\begin{align}
		\begin{split}
			\label{s1}
			&\delta\in \Pi \Leftrightarrow  \forall p,q\in\{-n,...,n\}, p\preceq q \text{ one has }:\\&\sum_{X\text{  contains }\geq (p+n+1-s(p,q))\text{ elements } \preceq q} \delta_{X}=m_{p,q},
		\end{split}
	\end{align}
	
	where  $s(p,q)$  was defined in   \eqref{snk}.

	\begin{lemma} There exists a one-to-one correspondence between the Gelfand-Tsetlin diagrams   $(m_{p,q})$  and shifted lattices $\Pi=\gamma+\mathcal{B}_{GC}^{\mathfrak{gl}_m}$ (using the system  \eqref{s1}) that contain at least one vector with only non-negative coordinates.
	\end{lemma}	
	
	The elements of the Gelfand-Tsetlin diagram, understood in the usual sense, are restored using the equalities \eqref{s1}.
	In this approach, diagrams related to the representation of the highest weight $[m_{-n},...,m_{n}]$ are defined as a class $mod\mathcal{B}_{GC}^{\mathfrak{gl}_m}$ of the vectors $\gamma$, such that  $\sum_{X: |X|=p}\gamma_X=m_{p+n-1-s(p)}$, $p\in\{-n,...,n\}$.
	
	Thus, {\it  it is possible to define } the Gelfand-Tsetlin diagram as a shifted lattice $\Pi=\gamma+\mathcal{B}_{GC}^{\mathfrak{gl}_m}$ of integer vectors, such that in this shifted lattice there is a representative having only non-negative coordinates.

	Then the Gelfand-Tsetlin diagram can be called $\mathfrak{gl}_{m-q}$-maximal if in its representative $\gamma$ (and, therefore, in all its vectors) the coordinates corresponding to subsets of $X$ turn to be zero, in the case when the monomial $A_X$ is not $\mathfrak{gl}_{m-q}$ is the highest vector under the action of  $\mathfrak{gl}_{m}$ on these variables defined in \eqref{dax}.

	\subsection{The GKZ base in the Zhelobenko's realization}
	
	In \cite{a4} the following was proved.

	\begin{theorem}
		\label{osa0}
		Consider $\Gamma$-series $\mathcal{F}_{\gamma}(A;\mathcal{B}_{GC}^{\mathfrak{gl}_m})$, such that $\Pi=\gamma+\mathcal{B}_{GC}^{\mathfrak{gl}_m}$ are Gelfand-Tsetlin diagrams in the sense of the previous section. Then $\mathcal{F}_{\gamma}(A;\mathcal{B}_{GC}^{\mathfrak{gl}_m})$ is the basis of the Zhelobenko's model $Zh$, and this basis is consistent in the decomposition of $Zh$ into the sum of irreducible representations.

	\end{theorem}

	\subsection{An A-GKZ model of representations  $\mathfrak{gl}_m$}
	
	Consider the action of the Lie algebra $\mathfrak{gl}_m$ on the independent variables $A_{X}$, defined by the formula
	
	\begin{equation}
		\label{dax}
		E_{i,j}=\sum_{X}A_{X,i}\frac{\partial}{\partial A_{X,j}}
	\end{equation}

	One has.
	
	\begin{theorem}[\cite{a4}]
		\label{osa}
		The space of polynomial solutions of the A-GKZ system is invariant under the action of $\mathfrak{gl}_m$. This space is a direct sum of finite-dimensional irreducible representations each taken with multiplicity $1$.
	\end{theorem}
	
	Thus, the solution space A-GKZ in the case of $\mathfrak{gl}_m$ is a model of finite-dimensional irreducible representations. It is not difficult to remark that this model is naturally identified with a subspace in the tensor product of standard representations.

	\subsection{A relation with the Gelfand-Tselin base }
	
	The following statement is proved in \cite{a4}. Introduce on the set of vectors considered by $mod\mathcal{B}_{GC}^{\mathfrak{gl}_m}$, the ordering
	
	\begin{equation}
		\label{poryadok}
		\gamma \preceq \delta \Leftrightarrow \gamma=\delta-s r\,\,\,\,\, mod  \mathcal{B}_{GC}^{\mathfrak{gl}_m},\,\,\,\, s\in \mathbb{Z}^{\mathcal{K}}_{\geq 0},
	\end{equation}
	
	where ${\mathcal{K}}$  is  the number of independent vectors  $v_{\alpha}$  of type  \eqref{va}, the vectors  $r_{\alpha}$  were defined in   \eqref{ra}. We use a notation  $sr:=s_1r_1+...+s_{\mathcal{K}}r_{\mathcal{K}}$.

	According to the section \ref{dgca}, the Gelfand-Tsetlin diagrams are in bijective correspondence with the shift vectors considered by $mod\mathcal{B}_{GC}^{\mathfrak{gl}_m}$, such that there is a vector with non-negative coordinates in this class. Then you can introduce a notation $G_{\delta}$ for the Gelfand-Tsetlin base vector corresponding to the diagram $\delta$.

	\begin{theorem}
		\begin{enumerate}

			\item The basis  $F_{\delta}$ and  $G_{\delta}$  are related by a lower-triangular transformation relatively the ordering   $\prec$.
			\item   $G_{\delta}$ is a orthogonalization of  $F_{\delta}$.
		\end{enumerate}
	\end{theorem}
	
	In  \cite{a4} a transformation that relates $F_{\delta}$ and  $G_{\delta}$ was found explicitly.

	\section{The GKZ and the A-GKZ systems for the series  $B$, $C$, $D$}
	
	\label{systemy}
	To obtain systems for which the analogs of the Theorems \ref{osa0}, \ref{osa} hold, one needs to add new equations to the GKZ and A-GKZ systems for the $A$ series. These equations are related to new (non-Plucker) relations arising for the series $B$, $C$, $D$.

	\subsection{ An auxiliary lattice  $\mathcal{B}^{g_n}$. Variables  $B_X$}

	Introduce a lattice
	
	\begin{equation}\bar{\mathcal{B}}^{g_n}=\mathcal{B}_{GC}^{\mathfrak{gl}_n}\subset \bar{L}= \mathbb{Z}^N,\,\,\,\,\, h_{\kappa}=e_X-e_{\widehat{-X}} ,\end{equation}
	
	The lattice $\bar{L}=\mathbb{Z}^N$ is interpreted as the set of exponents of Laurent monomials in variables $A_X$ introduced earlier.

	In the case of the series $C$ we define  an auxiliary lattice of $\mathcal{B}^{g_n}$ just as the lattice    $\bar{\mathcal{B}}^{g_n}$ defined above. But for the series $B$, $D$, further constructions are needed.

	Namely, consider variables $B_X$ associated with the previously introduced variables $A_X$ according to the following rule. If $X$ or $\widehat{-X}$ coincides with $\{\pm n,...,\pm 1\}$, then $B_{X}=\sqrt{A_{X}}$, in other cases $B_X=A_X$. When one passes to the Zhelobenko model in the first case $\sqrt{a_{X}}$ is substituted instead of $B _X$ (or $\sqrt{\bar{a}_{X}}$ for the $D$ series, when $|X|=n$, $m_{-1}<0$), and in the second case one substitutes just  $a_X$ (or $\bar{a}_{X}$ for the series $D$, for $|X|=n$, $m_{-1}<0$).
	Note that a construction of this type (a transition to variables, some of which coincide with the determinants, and some of which are the roots of the determinants in order for the spinor representation to be realized in the space of {\it polynomials} in
	these variables) for the case of $\mathfrak{o}_5$ is implemented in all details in \cite{a6}.

	Take the lattice $L=\mathbb{Z}^N$ of exponents of monomials in variables $B_X$.
	One can naturally embed the lattice  $\bar{L}$ in $L$.
	This gives the embedding $\bar{\mathcal{B}}^{g_n}\subset\bar{L}\subset L=\mathbb{Z}^N$.
	Define $\mathcal{B}^{g_n}$ as a sublattice in $L$, which is the image of $\bar{\mathcal{B}}^{g_n}$ for this embedding.

	Introduce the degrees of variables  $B_X$:
	\begin{equation}
		\label{podst}
		degB_X=1, \,\,\, X\neq \begin{cases}\{\pm n,...,\pm 1\},\\   \{\pm n,...,\pm 1,0\}  \end{cases},	degB_X=\frac{1}{2}, \,\,\,  \text{ otherwise }.
	\end{equation}
	
	\subsection{The GKZ system for  $g_n$ in the case of the series  $B$, $D$}
	
	Using the lattice $\mathcal{B}^{g_n}$, it is possible to build a system of GKZ  and then the A-GKZ system. But then one needs to add some new equations so that for the series $B$, $D$ the analogue of the statement formulated in the Remark \ref{rem1} holds. These new equations correspond to the Jacobi relations \eqref{sb}, as well as the root of the relations \eqref{s22}. Despite the addition of these equations, we  call the resulting systems {\it the GKZ and A-GKZ systems } for the series $B$, $D$.

	Consider functions in {\it independent} variables $B_X$, indexed by proper subsets in $\{-n,...,n\}$. In the case of the $B$ series, the index $0$  is included.

	\begin{definition}
		
		{\it The GKZ system for the series  $B$, $D$} is a system {\it generated } by  equations:
		
		\begin{equation}\label{ggkz}
			\begin{cases}
				& \mathcal{O}_{\alpha}\mathcal{F}:=(\frac{\partial^{\tau_1+\tau_2}}{\partial^{\tau_1} B_{i,Z}\partial^{\tau_2} B_{j,y,Z}}-\frac{\partial^{\tau_3+t_4}}{\partial^{\tau_3}B_{j,Z}\partial^{\tau_4} B_{i,y,Z}}\mathcal{F}=0,\\
				&  \tau_{i}=2,\text{ if the corresponding variable is of type $B_{\pm n,...,\pm 1}$ or  $B_{\pm n,...,\pm 1,0}$},\\& \text{ and $\tau_{i}=1$ otherwise},\\
				&(\frac{\partial}{\partial B_X}-\pm\frac{\partial}{\partial B_{\widehat{-X}}})\mathcal{F}=0,\text{  the sign is defined in   \eqref{znk}},\\&
				(\frac{\partial }{\partial B_{\pm n,...,\hat{\pm i},...,\hat{\pm j},...,\pm 1,-i,i}}-\frac{\partial^2 }{\partial   B_{\pm n,...,\hat{\pm i},...,\hat{\pm j},...,\pm 1,-i,-j} 
					\partial  B_{\pm n,...,\hat{\pm i},...,\hat{\pm j},...,\pm 1,i,j} })F=0,\\
				& (\frac{\partial^2}{\partial B_{\pm n,...,\hat{i},...,\pm 1 ,-i}\partial B_{\pm n,...,\hat{i},...,\pm 1 ,i}}-\frac{\sqrt{-1}}{\sqrt{2}}\frac{\partial}{\partial  B_{\pm n,...,\hat{i},...,\pm 1 ,0}})\mathcal{F}=0\text{ only for the series }B
			\end{cases}
		\end{equation}
	\end{definition}
	
	The equation of the second type is called  {\it the Jacobi equation}.

	\subsection{The A-GKZ system for the series  $B$, $D$}
	
	\label{ser1}

	\begin{definition}   {\it The A-GKZ system for the series  $B$, $D$} is a system {\it generated } by  equations:
		\begin{equation}
			\label{agkzgn}
			\begin{cases}
				&\bar{ \mathcal{O}}_{\alpha}\mathcal{F}:=(\frac{\partial^{\tau_1+\tau_2}}{\partial^{\tau_1} B_{i,Z}\partial^{\tau_2} B_{j,y,Z}}-\frac{\partial^{\tau_3+t_4}}{\partial^{\tau_3}B_{j,Z}\partial^{\tau_4} B_{i,y,Z}}+ \frac{\partial^{\tau_5+\tau_6}}{\partial^{\tau_5}B_{y,Z}\partial^{\tau_6} B_{i,j,Z}})F=0,\\
				&  \tau_i=2,\text{ if the corresponding variable is of type $B_{\pm n,...,\pm 1}$  or $B_{\pm n,...,\pm 1,0}$},\\& \text{ and $\tau_i=1$ otherwise},\\
				&(\frac{\partial}{\partial B_X}-\pm\frac{\partial}{\partial B_{\widehat{-X}}})F=0,\text{ the sign is defined in  \eqref{znk}},\\
				\\&
				(\frac{\partial }{\partial B_{\pm n,...,\hat{\pm i},...,\hat{\pm j},...,\pm 1,-i,i}}-\frac{\partial^2 }{\partial   B_{\pm n,...,\hat{\pm i},...,\hat{\pm j},...,\pm 1,-i,-j} 
					\partial  B_{\pm n,...,\hat{\pm i},...,\hat{\pm j},...,\pm 1,i,j} })F=0,\\
				& (\frac{\partial^2}{\partial B_{\pm n,...,\hat{i},...,\pm 1 ,-i}\partial B_{\pm n,...,\hat{i},...,\pm 1 ,i}}-\frac{\sqrt{-1}}{\sqrt{2}}\frac{\partial^2}{\partial^2  B_{\pm n,...,\hat{i},...,\pm 1 ,0}})F=0\text{ only for the series }B.
			\end{cases}
		\end{equation}
	\end{definition}
	\begin{remark}
		\label{rmrk4}
	The idea of construction of these  A-GKZ systems is the following. It is constructed  starting from  {\it some  } relations between  $a_X$, $\sqrt{a_{\pm n,...,
				\pm 1}}$.  In these relations the determinants are replaced firstly by variables   $A_X$,   and then by variables  $B_X$. Secondly a substitution    $B_X\mapsto \frac{\partial}{\partial B_X}$ is done. 
		
	 The first equation are constructed using 	 {\it some }  Plucker relations, the chosen relations correspond to the  vectors м  \eqref{va}.  The equations of the second type are written using the Jacobi relations. The equations of the third and the fourth type are written using the relations  \eqref{s22}, \eqref{s333}.
		
	 The GKZ system is a simplification of the  A-GKZ system.
	\end{remark}
	
	Let us construct basis in the spaces of solutions of the GKZ and A-GKZ systems.
	
	For an arbitrary vector $\delta\in\mathbb{Z}^N$ one introduces functions analogous to the functions \eqref{fs}.

	\begin{equation}
		\label{fss}
		\mathfrak{f}_{\delta}^s(B)=\sum_{t\in\mathbb{Z}^{\mathcal{K}}}\frac{(t+1)...(t+s)B^{\delta+tv}}{s!(\delta+tv)!},
	\end{equation}
	
	we use a chosen basis  $v_{\alpha}$ of the auxiliary lattice   $\mathcal{B}^{g_n}$.
	Put
	
	$$
	\mathfrak{f}_{\delta}(B):= \mathfrak{f}_{\delta}^0(B)
	$$
	
	Then
	
	\begin{equation}
		\label{reshenije}
		f_{\delta}(B):=\sum_{s\in \mathbb{Z}^{\mathcal{K}}_{\geq 0}}(-1)^s \mathfrak{f}_{\delta}^s(B)
	\end{equation}

	Both functions  $\mathfrak{f}_{\delta}(B)$  and $f_{\delta}(B)$ for  $\delta\in\mathbb{Z}^N$ are polynomials.
	
	The following ruler of differentiation holds
	
	\begin{equation}
		\label{prdf}
		\frac{\partial}{\partial B_X} \mathfrak{f}^s_{\delta}=\mathfrak{f}^s_{\delta-e_X},\,\,\,\frac{\partial}{\partial B_X} f_{\delta}=f_{\delta-e_X},
	\end{equation}

	where $e_X$ is the unit basis vector corresponding to the coordinate $B_X$.
	
	Also, as in the case of the $A$ series, it is proved that $\mathfrak{f}_{\delta}(B)$ and $f_{\delta}(B)$ are solutions of systems consisting of equations of the first type for the systems \eqref{ggkz} and \eqref{agkzgn} respectively.

	Introduce vectors

	\begin{align}
		\begin{split}
			\label{hkappa} 
			&h_{\kappa}=e_X-e_{\widehat{-X}},\,\,\,\,\,  \kappa=1,...,T\\
			&x_{\chi}=e_{\pm n,...,\hat{\pm i},...,\hat{\pm j},...,\pm 1,-i,i}-
			e_{\pm n,...,\hat{\pm i},...,\hat{\pm j},...,\pm 1,-i,-j}-\\&-e_{\pm n,...,\hat{\pm i},...,\hat{\pm j},...,\pm 1,i,j}, \,\,\,\,\,  \chi=1,...,Z
			\\&w_{\epsilon}=e_{\pm n,...,\hat{i},...,\pm 1 ,-i}+e_{\pm n,...,\hat{i},...,\pm 1 ,i}-e_{\pm n,...,\hat{i},...,\pm 1 ,0},\,\,\,\,\,\epsilon=1,...,E.
		\end{split}
	\end{align}

	Note that when one fixes $\kappa$,  then one fixes a subset of $X=\{i_1,...,i_t\}$. For a given $\kappa$ let the sign $(\pm_{\kappa} 1)$  be defined according to the ruler \eqref{znk}.	
	 Introduce the functions

	\begin{align}
		\begin{split}
			\label{sol1}
			&\mathcal{F}^s_{\delta}(B):=\sum_{t_1\in \mathbb{Z}^T,t_2\in\mathbb{Z}^{E},t_{3}\in \mathbb{Z}^{Z}} \Big (  \prod_{\kappa}(\pm_{\kappa} 1)^{t_1^{\kappa}}\cdot (\frac{\sqrt{-1}}{2})^{\sum_{\epsilon} t_2^{\epsilon}} \Big) \cdot \mathfrak{f}^s_{\delta+t_1 h+t_2 w+t_3x}(B ) \text{ for the series }B,\\
			&\mathcal{F}^s_{\delta}(B):=\sum_{t_1\in \mathbb{Z}^K,t_{3}\in \mathbb{Z}^{Z}} \Big (  \prod_{\kappa}(\pm_{\kappa} 1)^{t_1^{\kappa}} \Big) \cdot \mathfrak{f}^s_{\delta+t_1 h+t_3x}(B ) \text{  for the series }D
		\end{split}
	\end{align}
	\begin{equation}
		\label{sol}
		F_{\delta}(B):=\sum_{s\in\mathbb{Z}^{\mathcal{K}}_{\geq 0}}(-1)^s \mathcal{F}^s_{\delta}(B)
	\end{equation}

	One uses a notation $t_1 h:=t_1^{1} h_1+...+t_1^{T} h_T$,  $t_2 w:=t_2^{1} w_1+...+t_2^{E} w_E$, $t_3 \chi:=t_3^{1} x_1+...+t_3^{Z} x_{Z}$.
	Introduce a definition
	
	\begin{definition}
		The Gelfand-Tsetlin lattice is defined as follows
		\begin{align*}
			&\mathcal{B}^{g_n}_{GC}:=<\mathcal{B}^{g_n},h_{\kappa}, w_{\epsilon},x_{\chi}>\,\,\,\text{ for the series }B,\\
			&\mathcal{B}^{g_n}_{GC}:=<\mathcal{B}^{g_n},h_{\kappa},x_{\chi}>\,\,\,\text{ for the series }D
		\end{align*}
	\end{definition}

	Using  \eqref{prdf} one gets

	\begin{theorem}	
		\label{tagkz}
		The space of polynomial solutions of the GKZ system for the series $B$, $D$ has a basis consisting of  function	$\mathcal{F}_{\delta}(B)=\mathcal{F}^0_{\delta}(B)$
		for various $mod\mathcal{B}^{g_n}_{GC}$ vectors $\delta\in\mathbb{Z}^N$, such that this function is nonzero. 
		
		The space of polynomial solutions of the A-GKZ system for $g_n$ has the basis (with the same choice $\delta$) of the function $F_{\delta}(B)$

	\end{theorem}

	\begin{corollary}
		\label{gkzagkz}
		There is a one-to-one correspondence between the spaces of polynomial solutions of GKZ system \eqref{ggkz} and the A-GKZ system \eqref{agkzgn}.
	\end{corollary}	
	\begin{remark}
		This theorem is a manifestation of the following general idea in the analytical theory of differential equations: given a  system of equations, a simplified system is constructed, which is solved explicitly (GKZ is a "simplification" of A-GKZ). Then, for each solution of the simplified system, the solution of the initial system is constructed. There are several formalization of this idea  \cite{br}, \cite{st}.
	\end{remark}

	\begin{remark}
		A relation between the GKZ and A-GKZ systems is an example of the toric degeneration of differential operators.
		
		In the theory  between  of  the toric degeneration of Plucker relations in variables $A_X$ are considered. In this case, the basic relations are selected and some terms are removed  \cite{tg1}. This theory is related to  the representation theory,  the Gelfand-Tsetlin diagrams arise in it (see the original paper \cite{tg2} and the closet to the present paper   \cite{tg3}).
		
		Nevertheless, this relation  between the theory of toric degenerations  and the present paper  cannot be formalized. In the present paper, degenerations of differential operators are considered, spaces of their solutions are investigated, Gelfand-Tsetlin diagrams are used to index basic solutions. All this has no analogues in the works on toric expressions
		
	\end{remark}

	\subsection{The systems for the tensor representations in the case of the series $B$, $C$, $D$}
	\label{s2}
	
	It is known that in the case of the $C$ series, any finite-dimensional representation
	\footnote{Thus, the GKZ and A-GKZ systems for the $C$ series are the systems constructed in this Section.} is realized as a sub-representation in the tensor power of the standard representation. In the case of series
	$B$, $D$ it is possible to set the problem of construction of models only for such representations. To build a model of tensor representations, we will use more simple  GKZ and A-GKZ systems than in the previous section.

	Consider the independent variables $A_X$, antisymmetric in $X$. When one passes to the Zhelobenko model one substitutes $a_X$ instead of $A_X$  for all $X$.

	Take the auxiliary lattice $\bar{\mathcal{B}}^{g_n}$ constructed above. We denote it now also by  $\mathcal{B}^{g_n}$. Take a GKZ system based on it and add to it   the Jacobi equation.

	\begin{definition}   {\it The GKZ system for the tensor representations} is  the system {\it generated} by  equations
		
		\begin{equation}\label{ggkzt}
			\begin{cases}
				& \mathcal{O}_{\alpha}\mathcal{F}:=(\frac{\partial^2}{\partial A_{i,Z}\partial A_{j,y,Z}}-\frac{\partial^2}{\partial A_{j,Z}\partial A_{i,y,Z}})\mathcal{F}=0,\\
				&(\frac{\partial}{\partial A_X}-\pm\frac{\partial}{\partial A_{\widehat{-X}}})\mathcal{F}=0,\text{  the sign is defined in   \eqref{znk}}
			\end{cases}
		\end{equation}
		
	\end{definition}
	
	\begin{definition}  {\it The A-GKZ system for the tensor representations} is  the system {\it generated} by  equations
		
		\begin{equation}
			\label{agkzgnt}
			\begin{cases}
				&\bar{\mathcal{O}}_{\alpha}F=(\frac{\partial^2}{\partial A_{i,Z}\partial A_{j,y,Z}}-\frac{\partial^2}{\partial A_{j,Z}\partial A_{i,y,Z}}+\frac{\partial^2}{\partial A_{y,Z}\partial A_{i,j,Z}})F=0,\\
				&(\frac{\partial}{\partial A_X}-\pm\frac{\partial}{\partial A_{\widehat{-X}}})F=0.
			\end{cases}
		\end{equation}
	\end{definition}

	\begin{definition}
		The Gelfand-Tsetlin lattice for the tensor representations  is define as follows
		$$
		\mathcal{B}^{g_n}_{GC}:=<\mathcal{B}^{g_n},h_{\kappa}>
		$$
	\end{definition}
	
	The solutions of the A-GKZ system are constructed as follows 
	\begin{align}
		\begin{split}
			\label{sol10}
			&\mathcal{F}^s_{\delta}(A):=\sum_{t_1\in \mathbb{Z}^K} \Big (  \prod_{\kappa}(\pm_{\kappa} 1)^{t_1^{\kappa}} \Big) \cdot \mathfrak{f}^s_{\delta+t_1 h}(A ),
		\end{split}
	\end{align}
	
	where $\mathfrak{f}^s_{\delta}(A)$ is defined by the same formula  \eqref{fss},  where however one uses the variables  $A_X$ and another auxiliary lattice   $\mathcal{B}^{g_n}$.
	\begin{equation}
		\label{sol0}
		F_{\delta}(A):=\sum_{s\in\mathbb{Z}^{\mathcal{K}}_{\geq 0}}(-1)^s \mathcal{F}^s_{\delta}(A)
	\end{equation}
	
	An analogue of the Theorem \ref{tagkz} and the Corollary  \ref{gkzagkz} holds.

	\section{The Gelfand-Tsetlin diagrams for the algebra  $g_n$}
	\label{dgc}

	By analogy with the Section \ref{dgca} let use a definition of a Gelfand-Tsetlin diagram for the algebra $g_n$.
	
	One defines  Gelfand-Tsetlin diagrams as follows.

	\begin{definition}
		\label{d1}

		A Gelfand-Tsetlin diagram for the algebra $g_n$ is a shifted lattice $\Pi=\gamma+\mathcal{B}_{GC}^{g_n}$ consisting of integer vectors $\gamma\in\mathbb{Z}^N$, such that in $\Pi$ there is a vector that has only non-negative coordinates.

	\end{definition}
	
	In this case,  to a diagram $\delta$ there corresponds to a set of highest weights \begin{equation}\label{wtk}wt_{n-k+1}(\delta)=[m_{-n,k},...,m_{-1+(k-1),k}]\end{equation} for algebras $g_{n-k+1}$, $k=1,...,n$ constructed by the rule

	\begin{align}
		\begin{split}
			\label{strok}
			&m_{p,q}=\sum_{X: X\text{ contains } \geq (n+p+1)\text{ indices whose absolute value} \geq q}\gamma_X-\\&- \sum_{X: X\text{  contains } \geq (-p-1)\text{ indices whose absolute value  }  \geq q}\gamma_X,
			\\&p=-n,...,-1,\,\,\,\,\,\, q=1,...,n, p\leq -q+1
		\end{split}
	\end{align}
	
	This definition is correct, since it does not depend on the choice of the representative of $\gamma$ in the equivalence class $mod \mathcal{B}_{GC}^{\mathfrak{g}_n}$.

	\begin{definition}
		The row   $[m_{-n,1},..., m_{-1,1}]$ is called the highest weight of a diagram.
	\end{definition}
	
	Another definition of Gelfand-Tsetlin diagrams is more familiar, in which some tables are constructed that encode the basic vectors of a finite-dimensional irreducible representation of $g_n$ (see \cite{m}). At the same time, these tables contain the rows $m_{p,q}$, defined in \eqref{strok}.
	
	It is shown below that diagrams in the sense of Definition  \ref{d1} are also in one-to-one correspondence with the basic vectors of  finite-dimensional irreducible representations of $g_n$.
	From here one can conclude that the objects defined in \ref{d1} are manifestations of the same essence that is encoded by the Gelfand-Tsetlin diagram in the classical, combinatorial sense.
	
	The advantage of Definition  \ref{d1}  becomes clear when one calculates the formulas for the action of generators. If you think in terms of vectors $\gamma$, these formulas look very simple and natural (see  \eqref{deistviel}).

	\section{Models of representations in their cases  $B$, $C$, $D$}
	\label{mdl}
	
	\subsection{The action}
	
	\label{razddeist}
	
	There is an action of the algebra $g_n$ on the variables $A_X$, given by the formulas \eqref{deistviel}, where an expression $F_{i,j}$ through $E_{i,j}$ is used.
	
	The action of $g_n$ onto the variables $B_X$ in the case of the series $B$, $D$ for $X$ and $\widehat{-X}\neq \{\pm n,...,\pm 1\}$ is defined in the same way as onto $A_X$. In the case of $X$ or $\widehat{-X}=\{\pm n,...,\pm 1\}$, the action is defined so that $B_X$ for these $X$ form a spinor representation. To obtain such formulas, one needs to take the formulas of the action $g_n$ on the roots of the determinants $\sqrt{a_X}$ (or $\sqrt{\bar{a}_X}$. If one wants to obtain a representation with $m_{-1}<0$ for the series $D$) and replace $\sqrt{a_X}$ (or $\sqrt{\bar{a}_X}$) by $B_X$. Thus, in the case of the $D$ series, there are two ways to construct an action.

	The explicit form of the resulting formulas was discussed when proving the Sentence \ref{spb}, see formulas \eqref{fl1}, \eqref{chtn}.
	
	Thus, the spaces of polynomials $\mathbb{C}[A]$ and $\mathbb{C}[B]$ become representation spaces of the algebra $g_n$.

	\subsection{The plan}
	\label{plan}
	Next, we  get the following results.
	
	\begin{enumerate}
		\item  We prove that Gelfand-Tsetlin diagrams in the sense of Definition \ref{d1} encode the basic vectors of irreducible finite-dimensional representations.
		\item We prove that the solution space of the A-GKZ system is a  model of representations.
		\item  We construct a base of the Zhelobenko model consisting of $\Gamma$-series or monomials.
	\end{enumerate}

	To achieve these goals, we first prove that the solution space of the A-GKZ system coincides with the solution space of the ideal of relations between determinants. It follows from this that the solution space A-GKZ {\it contains } a representation model. Then we prove that the space of polynomial solutions of the A-GKZ model coincides with this model   (it is reasonable to call this model the A-GKZ model).  And thus the base solutions  $F_{\gamma}$  form a base of the A-GKZ model.

	As a corollary one gets that  the number of Gelfand-Tsetlin diagrams in the sense of Definition \ref{d1} for a fixed highest weight equals to the dimension of a representation with this highest weight.
	
	Then we pass to the  Zhelobenko's model. We prove that the span of $\Gamma$-series in determinants constructed from Gelfand-Tsetlin diagrams in the sense of \ref{d1} coincides with the Zhelobenko model. 
	From this we can already conclude that the number of Gelfand-Tsetlin diagrams in the sense of Definition \ref{d1} for a fixed highest weight {\it is equal to } the dimension of the representation.


	The goals are achieved in Section   \ref{cel2} (the first and the second goals); \ref{monomy} (the third goal ).
	
	\subsection{The A-GKZ model }
	\label{modela}
	\subsubsection{Invariance of the solution space of the A-GKZ system}

	This section discusses the A-GKZ systems defined in the Sections  \ref{ser1}, \ref{s2}.
	
Since the variables  $A_X$  are independent one can introduce an invariant scalar product  between polynomials  $f(A)$, $g(A)$ by the formula
	
	\begin{equation}
		\label{skp}
		<f(A),g(A)>:=f(\frac{\partial }{\partial A}) g(A)\mid_{A=0},
	\end{equation}

	where  $f(\frac{\partial }{\partial A}) $ is a result of substitution into    $f$  of differential operators  $\frac{\partial }{\partial A_X}$ instead of variables  $A_X$.
An analogous scalar product is introduced in the case of variables 	  $B$.

Now let us prove the following statement. 	
	
	\begin{propos}
		\label{prdl}
		The solution space of the A-GKZ system is $g_n$-invariant.
	\end{propos}
	
	\begin{proof}

	Let us consider separately the A-GKZ systems from Sections 	   \ref{ser1}, \ref{s2}.  in both cases  denote as   $\bar{I}$  the ideal that is contained in  $\mathbb{C}[\frac{\partial}{\partial B}]$   or   $ \mathbb{C}[\frac{\partial}{\partial A}]$  that corresponds to the considered A-GKZ system. 
		
		The action of   $g_n$ onto variables  $A_X$ induces an action on  $\mathbb{C}[\frac{\partial}{\partial A_X}]$. Analogously the action of  $g_n$ onto  $\mathbb{C}[\frac{\partial}{\partial B_X}]$ is defined.  
	Let us show that in both cases the ideal   $\bar{I}$ is invariant. 	 From this fact it immediately follows that the space $Sol_{\bar{I}_{}}$ of polynomial solutions of the considered  ideal is invariant.

	Before to proceed to consideration of the A-GKZ systems from Sections 	   \ref{ser1}, \ref{s2},  let us do an important notation. 
		
	Let us be given a system of partial differential equations corresponding to an ideal  $\bar{J}\subset \mathbb{C}[\frac{\partial}{\partial A}]$, which is obtained by a substitution  $A_X\mapsto \frac{\partial}{\partial A_X}$ from an ideal  $J\subset \mathbb{C}[A]$.   In terms of the scalar product    \eqref{skp}  the fact that $f(A)$ is a solution of a system of partial differential equations is equivalent to the fact that    $f(A)$ is orthogonal to the ideal $\bar{J}\subset \mathbb{C}[A]  $, corresponding to the system. Thus for an ideal $J\subset \mathbb{C}[A]$ one has $Sol_{\bar{J}}=(\bar{J})^{\perp}$. Note that   $(J^{\perp})^{\perp}$ is a closure of   $J$ in the topology induced by the scalar product   \eqref{skp}.   Due to the orthogonality of monomials one has $(J^{\perp})^{\perp}=J$.  Thus the coincidence of spaces of polynomial solutions for the system A-GKZ  for the series  $A$  and the ideal $\bar{I}_{\mathfrak{gl}_m}$  implies the coincidence of ideals 
		\begin{equation} \label{sovpad} \bar{I}_{\mathfrak{gl}_m}=\bar{I}^A   \subset \mathbb{C}[A]    ,\end{equation} 
		where $\bar{I}^A$  is an ideal that defines the A-GKZ system for the series    $A$.  From here one gets that the ideal   $\bar{I}^A$ is invariant.
		
		Let us prove the statement of the Proposition for the A-GKZ  system from Section  \ref{s2}. 
		Let  $I_{g_n}\subset \mathbb{C}[A]$  be an ideal of relations between determinants  $a_X$.  in Lemma	  \ref{lms} it was proved that the ideal    $I_{g_n}$  is generated by the ideal  $I_{\mathfrak{gl}_m}$ ($m=2n$ for the series  $C$, $D$,and  $m=2n+1$ for the series $B$)  and the Jacobi relations. 
		Hence,  $ \bar{I}_{g_n}$  is generated by the ideal $\bar{I}_{\mathfrak{gl}_m}$ and differential operators  that correspond to Jacobi relations.  Due to   \eqref{sovpad}  one gets that   $ \bar{I}_{g_n}$ is generated by  $\bar{I}^A$  and the operators that appear from the Jacobi relations. 
		But the ideal generated by    $\bar{I}^A$   and the Jacobi relations is just the ideal  $\bar{I}$ for the system
		\ref{s2}. That is for the system
		\ref{s2} one has
		
		\begin{equation}
			\label{sovpad2}
			\bar{I}=\bar{I}_{g_n}.
		\end{equation}

	As a consequence of   \eqref{sovpad2} the  ideal  $\bar{I}$ is invariant.

		Now consider the A-GKZ system for $g_n$ from the Section \ref{ser1}.
		Consider  a substitution of  variables $A_X$ instead of the variables $B_X$ by the rule:

		\begin{equation}
			\label{podst1}
			A_X\mapsto B_X, \,\,\, X\neq \begin{cases}\{\pm n,...,\pm 1\},\\   \{\pm n,...,\pm 1,0\}  \end{cases},	A_X\mapsto B^2_X, \,\,\,  \text{ othewise}.
		\end{equation}

		Directly (using transformations leading \eqref{ffa} to \eqref{fl1} and \eqref{chtn}), it is verified that this substitution is consistent with the action modulo the ideal $I'_{g_n}\subset \mathbb{C}[B]$ generated by the Plucker relations \eqref{spgl}, Jacobi relations and relations

		\begin{align}
			\begin{split}
				\label{sob}
				&B_{\pm n,...,\hat{i},...,\pm 1,i}B_{\pm n,...,\hat{i},...,\pm 1,-i}=\frac{\sqrt{-1}}{\sqrt{2}} B_{\pm n,...,\hat{i},...,\pm 1,i,0},\,\,\, \text{ only for the series  }B,\\
				& B_{\pm n,...,\hat{\pm i},...,\hat{\pm j},...,\pm 1,-i,i}=B_{\pm n,...,\hat{\pm i},...,\hat{\pm j},...,\pm 1,-i,-j}B_{\pm n,...,\hat{\pm i},...,\hat{\pm j},...,\pm 1,i,j}
			\end{split}
		\end{align}

		the square of which is a consequence of the Plucker and Jacobi relations.

		Now one can check that the ideal $\bar{I}'_{g_n}$ in the space $\mathbb{C}[\frac{\partial}{\partial B}]$ is invariant. For generators $\bar{I}'_{g_n}$ corresponding to the Plucker relations, this is checked using the equality \eqref{como} and the consistency of the substitution \eqref{podst} and the action $g_n$. For generators $\bar{I}'_{g_n}$ constructed by the Jacobi and \eqref{sob} relations, this is checked by direct calculation. Hence, the space of polynomial solutions $\bar{I'_{g_n}}$, which is nothing but the space of polynomial solutions of the A-GKZ system from the section \ref{s2}, is also invariant.

	\end{proof}

	\begin{lemma}\label{ostnach} The solution space of the A-GKZ system contains a model of representations.
		
	\end{lemma}
	
	\begin{proof}
		
		According to the Proposition \ref{prdl}, the solution space of the A-GKZ system for $g_n$ is a representation of $g_n$.
		Consider in the space $Sol_{AGKZ}$ solutions of A-GKZ for $g_n$ a finite-dimensional subspace of functions $Sol_{AGKZ}^{l_{-n},...,l_{-1}}$, of a fixed homogeneous degree $l_{p}$ by $A_X$ (or $B_X$) with $|X|=(p+n+1)$ and $(-p-1)$, $p=1,...,n$. One has

		$$
		Sol_{AGKZ}=\bigoplus_{l_{-n},...,l_{-1}}Sol_{AGKZ}^{l_{-n},...,l_{-1}}
		$$

		Under the action of $g_n$, homogeneous powers of determinants of the same size are preserved. Note also that the basis \eqref{sol} is consistent with this decomposition of the solution space A-GKZ for $g_n$

		Put $m_{p}=l_p+l_{p+1}+...$.
		Then in   $Sol_{AGKZ}^{l_{-n},...,l_{-1}}$  the following $g_n$-highest vector of the weight  $[m_{-n},m_{-n+1},...,m_{-1}]$ is contained
		
		\begin{equation}
			\label{stv}
			(A_{-n}+\pm A_{\widehat{n}})^{m_{-n}-m{-n+1}}( A_{-n,-n+1}+\pm A_{\widehat{n-1,n}})^{m_{-n+1}-m_{-n+2}}....,
		\end{equation}
		
		Or, in the case of series $B$, $D$ is a similar expression with variables $B$. Moreover, in the case of the $D$ series
		the last $n$-th factor is taken to the power  $|m_{-1}|$. The signs in brackets are defined in \eqref{znk}.
		
		Thus  $Sol_{AGKZ}^{l_{-n},...,l_{-1}}$ contains an irreducible representation of the highest weight $[m_{-n},m_{-n+1},...,m_{-1}]$.

	\end{proof}

	\subsubsection{The minimality of the solutioin space of the A-GKZ system}
	
	Let us prove the followings statement.
	
	\begin{lemma}  The solution space of the A-GKZ system and the  model of representations constructed in the proof of  Lemma  \ref{ostnach}    coincide.
	\end{lemma}
	Let us introduce for the model of representations constructed in the proof of   Lemma  \ref{ostnach}  the notation    $Mod$.
	
	\begin{proof}
		Let us suggest that 	  $Mod\neq Sol_{AGKZ}$.
		When one substitutes instead of independent variables 	$A_X$ (or   $B_X$, let us for simplicity below consider the case of variables $A_X$, the case when $f$ depends on $B_X$ is considered analogously) the determinants  $a_X$  the representation  $Sol_{AGKZ}$ is mapped to the Zhelobenko's  model and the subrepresentation   $Mod$ is mapped to the Zhelobenko's model isomorphically.  From here one gets that  (for example by considering the finite dimensional subspaces with the fixed homogeneous degree by  $A_X$, $|X|=l_i$) that the kernel of this mapping is non-trivial. Thus there exists non-zero function  $f\in Sol_{AGKZ}$  such that
		
		$$
		f\mid_{A_X \mapsto a_X}=0,
		$$
		
		that is  $f(A_X)\in I_{g_n}$.  But $f(\frac{\partial}{\partial A})g(A)=0$ for all  $g\in Sol_{AGKZ}$.  In particular  $f(\frac{\partial}{\partial A})f(A)=0$,  but this is possible only when   $f=0$.  Thus we obtain a contradiction.
	\end{proof}
	Since in 
	$Sol_{AGKZ}^{l_{-n},...,l_{-1}}$  there exists a base  \eqref{sol},  indexed by Gelfand-Tsetlin diagrams for  $g_n$  in the sence of Definition   \ref{d1}, then one comes to a conclusion

	\begin{corollary}
		\label{grsn}
		Fix a highest weight $[m_{-n},...,m_{-1}]$. Then the number of Gelfand-Tsetlin diagrams in the sense of \ref{d1} with such a highest weight i equals to the dimension of this irreducible representation.
	\end{corollary}

	\subsection{The A-GKZ model }
	\label{cel2}
	
	Thus the solution space of the A-GKZ system is a model of representations of  $g_n$,  the highest weight is given by the formulas \eqref{stv}.  Let us formulate this result as the following Theorem.

	\begin{theorem}\label{osnt0} 
		
		There exists an A-GKZ model formed by the space of solutions of the A-GKZ system. In the case of series $B$, $D$, it has a basis consisting  of functions $F_{\gamma}(B)$, and in the case of series $C$ it has a basis consisting of functions $F_{\gamma}(A)$.One takes as $\gamma$ all possible Gelfand-Tsetlin diagrams in the sense of the definition of \ref{d1}. Those functions for which
		the diagram has the highest weight $[m_{-n},...,m_{-1}]$,
	form a basis in the representation space of the highest weight $[m_{-n},...,m_{-1}]$.

	\end{theorem}
	




	Thus, the  first and the second goals from the \ref{plan} section are achieved.

	\subsection{ $\Gamma$-series in the Zhelobenko's realization}
	\label{zhlb}
	
	We use a short notation 
	
	$$
	\mathcal{F}_{\gamma}:=\mathcal{F}_{\gamma}(A) \text{ or }\mathcal{F}_{\gamma}(B),\,\,\,
	F_{\gamma}:=F_{\gamma}(A) \text{ or }F_{\gamma}(B).
	$$
	
	\subsubsection{A formulation of the result: $\Gamma$-series generate the Zhelobenko's realization}

	Consider $\Gamma$-series that are solutions of the GKZ system for the corresponding algebra $g_n$. Substitute the determinants $a_X$ in them instead of $B_X$ and $A_X$ according to the substitution rule described above. Let $W\subset Fun$ be their linear span.
	
	\begin{theorem}
	\label{gkzt}
	The space $W$ is a representation model of the algebra $g_n$.
	\end{theorem}
	From this Theorem, the coincidence of $W$ and the Zhelobenko model immediately follows.
	
	\begin{proof}
	Consider the functions $\mathcal{F}_{\gamma}$ with a  fixed homogeneous degree $l_i$ by determinants of the same size $i$. Let's check that when the generator $F_{i,j}$ acts on a function $\mathcal{F}_{\gamma}$, a linear combination of the functions $\mathcal{F}_{\delta}$ of the same homogeneous degree is obtained.
	To do this, we will use the Principal Lemma.

	\subsubsection{The formulation of the Principal Lemma}

	\begin{lemma} 
		\label{osnovl}	
		Let  us be given a $\Gamma$ series  $\mathfrak{f}_{\gamma}(c)$ in variables $c_X=a_{X}$ or $\sqrt{a_X}$. The lattice on which this series is constructed has a basis $v_{\alpha}$, $\alpha=1,...,K$, such that for $\alpha=1,...,{\mathcal{K}}$, ${\mathcal{K}}<K$ these generators have the following properties. They have the form

		$$
		v_{\alpha}=\tau_1\cdot e_{X_1}+\tau_2\cdot e_{X_2}-\tau_3\cdot e_{X_3}-\tau_4\cdot e_{X_4},
		$$
		where  $e_X$  is a unit vector corresponding to the coordinate $a_X$, $\tau_i\in\mathbb{Z}_{\geq 0}$.
		There exists a relation between variables of the form is associated with each such generative
		
		\begin{equation}
			\label{sbv}
			c^{\tau_1}_{X_1}c^{\tau_2}_{X_2}-c^{\tau_3}_{X_3}c^{\tau_4}_{X_4}+c^{\tau_5}_{X_5}c^{\tau_6}_{X_6}=0,
		\end{equation}

		for some  $X_5,X_6$, $\tau_5$, $\tau_6$. With each generator  $v_{\alpha}$  one assocates a vector
		
		$$
		r_{\alpha}=\tau_1\cdot e_{X_1}+\tau_2\cdot e_{X_2}-\tau_5\cdot e_{X_5}-\tau_6\cdot e_{X_6}.
		$$
		
		Then modulo the relations   \eqref{sbv} one has

		\begin{equation}
			\label{umn}
			c_X\mathfrak{f}_{\gamma}(c)=\sum_{s\in\mathbb{Z}^{\mathcal{K}}_{\geq 0}} C^{\gamma}_s\mathfrak{f}_{\gamma+e_X+sr}(c).
		\end{equation}
		
		
	\end{lemma}

	\subsubsection{The proof of the Principle Lemma }

	Take the independent variables $C_X$, antisymmetric by $X$.
	For a polynomial $g(C)$  denote as $$g(\frac{\partial}{\partial C})$$  a differential operator obtained by replacing each variable $C_{X}$ by differentiating $\frac{\partial}{\partial C_X}$.

	Denote as $Plk$ the ideal generated by the relations

	$$
	C^{\tau_1}_{X_1}C^{\tau_2}_{X_2}-C^{\tau_3}_{X_3}C^{\tau_4}_{X_4}+C^{\tau_5}_{X_5}C^{\tau_6}_{X_6}=0,
	$$

	With this relation, one relates a system A-GKZ, {\it  generated } by  equations.
	\begin{align*}
		\Big( (\frac{\partial }{\partial  C_{X_1}})^{\tau_1} (\frac{\partial }{\partial  C_{X_2}})^{\tau_2}-(\frac{\partial }{\partial  C_{X_3}})^{\tau_3} (\frac{\partial }{\partial  C_{X_4}})^{\tau_4}+(\frac{\partial }{\partial  C_{X_5}})^{\tau_5} (\frac{\partial }{\partial  C_{X_6}})^{\tau_6} \Big)f=0,
	\end{align*}
	
	and  the equation of the system of GKZ corresponding to the rest of the generating  $v_{\alpha}$, $\alpha=k+1,...,K$.

	There is a basis in the space of solution of  this system
	
	$$f_{\delta}(C)=\sum_{s\in \mathbb{Z}_{\geq 0}^{\mathcal{K}}}(-1)^s\mathfrak{f}_{\delta-sr}^s(C)$$ РІ where $\mathfrak{f}_{\delta-sr}^s(C)$ is defined by analogy with \eqref{fss}.

	Let us do the following observation (similar to how it is done in \cite{a4})

	\begin{lemma}
		$\lambda_1 g_1(c)+...+\lambda_l g_l(c)=0\,\,\, mod Plk $, if
		$$
		(\lambda_1 g_1(\frac{\partial }{ \partial C})+...+\lambda_l g_l(\frac{\partial }{\partial C}))f_{\delta}(C)=0
		$$
		for all  $C$.
	\end{lemma}

	Put  $\delta+v:=\delta+v_1+...+v_{\mathcal{K}}$.
	Literally, as well as in \cite{a4}, the following statement is proved
	
	\begin{lemma}
		\label{lma}
		\begin{equation}\label{rwo}
			\mathfrak{f}_{\gamma}(\frac{\partial }{\partial C})f_{\delta}(C)=\sum_{s\in\mathbb{Z}^{\mathcal{K}}_{\geq 0}}\mathfrak{f}_{\gamma+v+sr}^s(1)f_{\delta-\gamma-sr}(C),
		\end{equation}
		
		where  $\mathfrak{f}_{\gamma+v+sr}^s(1)$ is a   result of substitution of units instead of all arguments  $C$.
	\end{lemma}

	Let us show that a similar equality exists for the functions $\mathcal{F}_{\gamma}(c)$ defined in \eqref{sol1}. Indeed, this function is the sum of $\Gamma$-series of the form $\mathfrak{f}_{\delta}$, for each of which has equality of the form
	\eqref{rwo}. Summing them up, one gets the following equality
	
	\begin{equation}\label{rwo1}
		\mathcal{F}_{\gamma}(\frac{\partial }{\partial C})f_{\delta}(C)=\sum_{s\in\mathbb{Z}^{\mathcal{K}}_{\geq 0}}\mathcal{F}_{\gamma+v+sr}^s(1)f_{\delta-\gamma-sr}(C),
	\end{equation}
	
	From here, literally as well as in \cite{a4}, we can conclude
	
	\begin{equation}
		\label{bxf}
		c_X\mathcal{F}_{\gamma}(c)=\sum_{s\in\mathbb{Z}^k_{\geq 0}} const^{\gamma}_s\cdot \mathcal{F}_{\gamma+e_X+sr}(c)
	\end{equation}
	
	The Principal Lemma is proved.

	\subsubsection{Coefficient in the Principal Lemma}
	
	One can give a formula for the coefficients in the expansion \eqref{umn}. Literally, as in \cite{a4}, one gets:

	\begin{align}
		\begin{split}
			\label{cs}
			&const^{\gamma}_s=\frac{\mathcal{F}_{\gamma+v}^{s}(1)}{\mathcal{F}^s_{\gamma+v+e_X+sr}(1)}-\sum_{p\in \mathbb{Z}^{\mathcal{K}}_{\geq 0}: s-p\in \mathbb{Z}^{\mathcal{K}}_{\geq 0},p\neq s}\frac{\mathcal{F}_{\gamma+v}^{p}(1)
				\mathcal{F}_{\gamma+v+pr+e_X+(s-p)r}^{s-p}(1)}{\mathcal{F}_{\gamma+v+pr+e_X+(s-p)r}(1)\mathcal{F}_{\gamma+v+pr+e_X}(1)}=\\
			&=
			\frac{\mathcal{F}_{\gamma+v}^{s}(1)}{\mathcal{F}^s_{\gamma+v+e_X+sr} (1) }-\sum_{p\in \mathbb{Z}^{\mathcal{K}}_{\geq 0}: s-p\in \mathbb{Z}^{\mathcal{K}}_{\geq 0},p\neq s}\frac{\mathcal{F}_{\gamma+v}^{p}(1)\mathcal{F}_{\gamma+v+e_X+sr}^{s-p}(1)}{\mathcal{F}_{\gamma+v+e_X+sr}(1)
				\mathcal{F}_{\gamma+v+pr+e_X}(1)}
		\end{split}
	\end{align}

	\subsubsection{The action of generators. The end of the proof of  Theorem \ref{gkzt}}
	
	Recall that $c_X=a_X$ or $\sqrt{a_X}$. The action of $F_{i,j}$ on $a_X$ is given by the formula \eqref{edet1}, the action on $\sqrt{a_X}$ is described in the proof of Proposition \ref{spb}, see formulas \eqref{fl1}.
	
	In both cases one can write

	\begin{equation}
		F_{i,j}=\sum_{Y_1,Y_2} c_{Y_1}\frac{\partial}{\partial c_{Y_2}}.
	\end{equation}

	Using the rule of differentiation  of  $\Gamma$-series and the Principle Lemma \ref{osnovl}, one   gets
	\begin{lemma}
		
		\begin{equation}
			\label{deistviel}
			F_{i,j}\mathcal{F}_{\gamma}=\sum_{Y_1}\sum_{s\in \mathbb{Z}^{\mathcal{K}}_{\geq\ 0}}
			const_{s}^{\gamma-e_{Y_2}} \cdot \mathcal{F}_{\gamma-e_{Y_2}+e_{Y_1}+sr}(c)
		\end{equation}

	\end{lemma}
	
	\begin{corollary}
		The linear span of functions $\mathcal{F}_{\gamma}$ is a representation.
	\end{corollary}
	
	Let's return to the proof of the Theorem \ref{gkzt} . In the Theorem \ref{ost1}, the conditions for the function necessary and sufficient for it to be included in the Zhelobenko model are written out.

	Applying it, one sees that $W\subset Yar$ is a sub-representation containing all finite-dimensional irreducible representations. Then $W=Zh$.
	The Theorem \ref{gkzt} is proved.

	\end{proof}

	As a corollary one gets.
	
	\begin{corollary}
	\label{grsv}
	
	Fix the highest weight of $[m_{-n},...,m_{-1}]$. Then the number of Gelfand-Tsetlin diagrams in the sense of Definition \ref{d1} with this highest weight is not less than the dimension of this irreducible representation.
	
	\end{corollary}

	\subsubsection{The Gelfand-Tsetlin diagrams} 
	
	\label{cel1}
	
	Using the Corollaries  \ref{grsn}, \ref{grsv}, one gets.
	
	\begin{theorem}\label{chdrzm}
	
	Fix the highest weight  $[m_{-n},...,m_{-1}]$. Then the number of Gelfand-Tsetlin diagrams in the sense of Definition \ref{d1} with this weight is equal to the dimension of this irreducible representation.
	
	\end{theorem}
	
	Thus, the first goal from the \ref{plan} section has been achieved.

	\subsubsection{Basis in the Zhelobenko's realization}
	\label{monomy}
	
	In the previous Sections, the functions $\mathcal{F}_{\gamma}(A)$ and $\mathcal{F}_{\gamma}(B)$ were constructed. In Sections \ref{ser1}, \ref{s2} when defining variables $A_X$, $B_X$ it was also indicated which substitution of determinants in these variables should be carried out when  one passes to the Zhelobenko's realization. Let us do such a substitution. The resulting function in both cases is denoted as $\mathcal{F}_{\gamma}(a) $.

	\begin{theorem}
	
	The functions $\mathcal{F}_{\gamma}(a)$, constructed for various Gelfand-Tsetlin diagrams in the sense of the Definition  \ref{d1}, form the basis of the Zhelobenko's realization.

	\end{theorem}
	
	\begin{proof}
	
	Since all these functions satisfy the conditions of the Theorem  \ref{gkzt}, they all belong to the Zhelobenko's model and it remains to prove the linear independence of these functions.
	
	The linear span of the function $\mathcal{F}_{\gamma}(a)$, for Gelfand-Tsetlin diagrams in the sense of  Definition  \ref{d1}, having the highest weight $[m_{-n},...,m_{-1}]$, form a sub-representation in the Zhelobenko's model of the highest weight $[m_{-n},...,m_{-1}]$. Their number is equal to the dimension of this representation. So they are independent, and therefore form a basis.

	\end{proof}

	Note that we not only proved that the linear span of $\mathcal{F}_{\gamma}(a)$ is a  model of representations of $g_n$, but also we explicitly wrote out formulas for the action of algebra generators.

	Let's prove another  statement. When defining the Gelfand-Tsetlin lattice $\mathcal{B}_{GC}^{g_n}$, its embedding into the lattice $\mathbb{Z}^N$ of exponents (integers or half-integers) of monomials in determinants was given (see Section \ref{systemy}). The image of the vector $\gamma\in \mathcal{B}_{GC}^{g_n} $ with this embedding is denoted as $\bar{\gamma}$. A similar construction can be implemented with vectors of shifted lattices. So for a Gelfand-Tsetlin diagram $\Pi$, we choose a representative of $\gamma$, then we construct a monomial in  determinants $a^{\bar{\gamma}}$.

	\begin{theorem}
	
	For every Gelfand-Tsetlin diagram in the sense of Definition  \ref{d1}, we construct a monomial $a^{\bar{\gamma}}$ according to the above procedure. Then these monomials form a basis in the Zhelobenko's realization.

	\end{theorem}
	\begin{proof}
	
	Indeed, by Lemma \ref{lma} one has
	
	$$
	\mathcal{F}_{\gamma}(a)=a^{\bar{\gamma}}+h.o.t.,
	$$
	where  $h.o.t.$   is a sum of monomials  $a^{\bar{\delta}}$, and  $\gamma\prec \delta$, where the ordering $\prec$ defined in the same way as in \eqref{poryadok}, with the replacement of the lattice by $\mathcal{B}_{GC}^{g_n}$.
	
	In particular one sees that  if $\gamma=\gamma'\,\,\, mod \mathcal{B}_{GC}^{g_n}$, then
	
	$$
	a^{\bar{\gamma}}=a^{\bar{\gamma'}}+h.o.t.
	$$
	
	Consider a set of vectors $\gamma\in\mathbb{Z}_{\geq 0}^N$, which are representatives of Gelfand-Tsetlin diagrams with the highest weight $[m_{-n},...,m_{-1}]$. Then the functions $\mathcal{F}_{\gamma}(a)$ form a basis
	of this representation in the Zhelobenko's model, and the set of functions $a^{\bar{\gamma}}$ is related with $\mathcal{F}_{\gamma}(a)$ by an upper-triangular transformation. Hence, the set $a^{\bar{\gamma}}$ also forms a basis.

	\end{proof}
	Thus, the third goal from the \ref{plan} section has been achieved.

	\section{ A relation to the Gelfand-Tsetlin base}
	\label{gttl}
	Let's establish a relation between the bases $F_{\gamma}(A)$ (or $F_{\gamma}(B)$) in the A-GKZ realization and the basis $\mathcal{F}_{\gamma}(a)$ in the Zhelobenko's realization with the Gelfand-Tsetlin basis $G_{\gamma}$.

	A Gelfand-Tsetlin type base can be defined as an eigenbasiss for the following maximal commutative subalgebra in $GT\subset U(g_n)$, called the Gelfand-Tsetlin algebra. For the subalgebra chain $g_1\subset g_2\subset ...\subset g_n$ take the centers of universal wrappers $Z(U(g_1)),...,Z(U(g_n))$ and generate a subalgebra with them in $GT\subset U(g_n)$ This is the Gelfand-Tsetlin algebra. You can explicitly write out the generators

	\begin{equation}
	\label{sso}
	C_{p}^q=\sum_{i_1,...,i_p\text{ modulo }\geq q}F_{i_1,i_2}F_{i_2,i_3}...F_{i_{2p},i_1}
	\end{equation}

	It is directly verified that these generators are self-adjoint with respect to the scalar product \eqref{skp}.
	
	It is known that the eigen-basis for the algebra $GT$, that is, the Gelfand-Tsetlin basis, is not unique in the case of the series $B$, $C$, $D$. Let us  give a construction of some Gelfand-Tsetlin basis. To do this, we first prove the following statement

	\begin{lemma}
	\label{lmff}
	The scalar product $<\mathcal{F}_{\gamma_1}(a),\mathcal{F}_{\gamma_2}(a)>$ can be nonzero only if there exists $\omega$ such that $\gamma_1\preceq \omega$, $\gamma_2\preceq\omega$.
	
	\end{lemma}
	
	\begin{proof}

	For certainty, we  conduct arguments using variables $A_X$. When using variables $B_X$, the reasoning is literally the same.
	
	To calculate the scalar product, let's go to the A-GKZ realization. Then the vector $\mathcal{F}_{\gamma_1}(a)$ in the Zhelobenko's realization  is written in A-GKZ  realization as a function of the form $\mathcal{FF}_{\gamma_1}(A):=\mathcal{F}_{\gamma_1}(A)+h(A)$, where $h(A)\in I_{g_n}$. Using the formula \eqref{skp} and taking into account that $F_{\delta}(A)$ is vanished by the ideal  $\bar{I}_{g_n}$, one gets

	$$
	<\mathcal{F}_{\gamma_1}(A)+h(A),F_{\delta}(A)>=	<\mathcal{F}_{\gamma_1}(A),F_{\delta}(A)>.
	$$
	
	By definition  $<\mathcal{F}_{\gamma_1}(A),F_{\delta}(A)>$ can be non-zero only if   $\gamma_1\preceq \delta$. Then  $\mathcal{FF}_{\gamma_1}=\sum_{s\in \mathbb{Z}^{\mathcal{K}}_{\geq 0}} c_{\gamma_1}^s F_{\gamma_1+sr}(A)$. Thus
	$$<\mathcal{F}_{\gamma_1}(a),\mathcal{F}_{\gamma_2}(a)>=\sum_{s_1,s_1\in\mathbb{Z}^{\mathcal{K}}_{\geq 0}   }c_{\gamma_1}^{s_1} c_{\gamma_2}^{s_2} <F_{\gamma_1+s_1r}(A),  F_{\gamma_2+s_2r}(A)  > .$$

	Considering  the supports, one can  conclude that the scalar product of $F_{\gamma_1+s_1r}(A)$ and $F_{\gamma_2+s_2r}(A) $ can be nonzero only if there exists $\omega$ such that $\gamma_1+s_1r, \gamma_2+s_2r\preceq\omega$
	This condition is equivalent to the condition from the formulation of the Lemma.

	\end{proof}

	\begin{corollary}
	
	There exists an orthogonal basis $\mathcal{G}_{\gamma}$ in the Zhelobenko's model, expressed in terms of the basis $\mathcal{F}_{\gamma}$ in an upper-triangular way with respect to the ordering $\prec$, that is, $$\mathcal{G}_{\gamma}(a)=\sum_{s\in\mathbb{Z}^{\mathcal{K}}_{\geq 0}} d_{\gamma}^s\cdot\mathcal{F}_{\gamma+sr}(a),\,\,\,\, d_{\gamma}^0=1.$$
	
	\end{corollary}

	\begin{theorem}
	\label{tgc}
	The base $\mathcal{G}_{\gamma}$ is the Gelfand-Tsetlin base.
	\end{theorem}
	\begin{proof}

	The Gelfand-Tsetlin base is an eigenbase for the Gelfand-
	Tsetlin algebra. Its generators \eqref{sso} are self-adjoint with respect to the scalar product \eqref{skp}, so that the space of a finite-dimensional irreducible representation $V\subset Zh$ is represented as an orthogonal direct sum of eigenspaces for $GT$. The set of eigenvalues defining one of the direct summands is given by the set of $g_{n-k}$-higher weights arising from the decomposition of $V$ into the sum of irreducible representation under the standard procedure of  restriction of algebras $g_n\downarrow g_{n-k}$. Introduce a notation $mu=\{\mu_n,...,\mu_1\}$ for the set of $g_n,...,g_1$-higher weights. The corresponding term in $V$ is denoted as $V_{\mu}$.

	Also, for the Gelfand-Tsetlin diagram $\gamma$, we introduce the notation (see \eqref{wtk})

	$$
	V^{\gamma}:=V_{\mu},\,\,\,\text{ where } \mu=\{wt_{n}(\delta),....,wt_{1}(\delta)\}.
	$$
	
	Let us prove the following statement.
	
	\begin{propos}
		\label{pfv}
		$\mathcal{F}_{\gamma}(a)\in\oplus_{s\in\mathbb{Z}^{\mathcal{K}}_{\geq 0}} V^{\gamma+sr}$
	\end{propos}
	
	\begin{proof}

		Consider first the monomial $a^{\delta}$. Applying the Jacobi relation if necessary, we can assume that $a^{\delta}$ depends only on the determinants $a_X$ with $|X|\leq n$.

		Let's introduce the {\it raising} operation.
		Let us first consider a realization in the space of functions on the subgroup $Z$ of upper-unitriangular matrices on the corresponding group $G$ (see \cite{zh}). Denote as $z_{i,j}$  functions of matrix elements on $Z$. In this case, a finite-dimensional irreducible representation of $V$ is realized in the space of polynomials in variables $z_{i,j}$ (including in the case of a half-integer of the highest weight) satisfying the indicator system \eqref{indsys}, the exponents in which are determined by the highest weight according to the rule \eqref{rb}.

		The monomial $a^{\delta}$ being restricted to $Z$ is written as a function $f(z_{i,j})$. In this case, one can assume that
		$i<-j$ in the case of series $B$, $D$ and $i\leq -j$ in the case of series $C$. The remaining variables $z_{i,j}$ are expressed in terms of these variables as polynomials.
		
		A function that is a $g_{n-k}$-the highest vector depends only on the variables $z_{i,j}$, $j\in\{-k,...,\hat{0},...,k\}$.


		Let us define the {\it raising} procedure as follows.  One applies the  operators $F_{i,j}$, $i<j$ until one gets a $g_{n-k}$-highest vector. Let us agree that this is done in the following order. First one applies $F_{-n,-(n-1)}$ to the maximum  possible power  (until the result is a non-zero), then $F_{-n,-(n-2)}$ to the maximum  possible power, then $F_{-(n-1),-(n-2)}$, etc. That is, the operators are applied according to the order
		
		\begin{equation}
			\label{proprem}
			F_{-n,-(n-1)};  F_{-n,-(n-2)} ; F_{-(n-1),-(n-2)} ; F_{-(n,-(n-3))};...
		\end{equation}

		In general the positive root operators $F_{i,j}$, $i<j$ in the realization under consideration are written as follows

		$$
		F_{i,j}=\frac{\partial}{\partial z_{i,j}}+\sum_{t<i}z_{t,i}\frac{\partial }{\partial z_{t,j}}-\pm\frac{\partial}{\partial z_{-j,-i}}+\sum_{t<-j}z_{t,-i}\frac{\partial }{\partial z_{t,j}},
		$$	
		
		where $\pm=+$ for series $B$, $D$ and $sign(i)sign(j)$ for series $C$. But since we apply   the operators $F_{i,j}$ in a certain order they are written simply as

		$$
		F_{i,j}=\frac{\partial}{\partial z_{i,j}}-\pm\frac{\partial}{\partial z_{-j,-i}}.
		$$

		And taking into account on which $z_{i,j}$ our function depends on, one has $F_{i,j}=\frac{\partial}{\partial z_{i,j}}$.
		
		The action of a {\it raising } procedure  onto a monomial $a^{\delta}$ can be described without going to the realization in the space of functions on $Z$. This operation acts as a substitution for $a_X\mapsto a_{X'}$, where $X'$ is obtained by the maximum possible left shift of all indexes whose absolute value $\geq k$.

		Each function $f(z_{i,j})$ can be written as follows
		
\begin{equation}
	\label{fpovysh}
	f=\sum_{\beta}c_{\beta}\cdot \prod_{|s|>k}\frac{z_{r,s}^{\beta_{r,s}}}{\beta_{r,s}!}\cdot f_{\beta}(z_{i,j}),
\end{equation}

		where $\beta$ is some index listing the terms in $f$ of the specified type.
		In this case, the function $f_{\beta}(z_{i,j})$ depends only on the variables $z_{i,j}$, $j\in \{-k,...,\hat{0},...,k\}$.
		When applying the {\it raising } procedure  in the expression \eqref{fpovysh} one gets $\sum_{\beta}c_{\beta}f_{\beta}(z_{i,j})$, where the sum of $\beta$ is taken, such that the vector

		$$
		[\alpha]:=(\alpha_{-n,-(n-1)}+\alpha_{(n-1),n},\alpha_{-n,-(n-2)}+\alpha_{(n-2),n}, \alpha_{-(n-1),-(n-2)}+\alpha_{(n-2),(n-1)},...)
		$$

		is the maximum relative to the lexicographic order. Such terms we call the  {\it maximal}. According to our assumption, on which $z_{i,j}$ the function $f$ depends, made at the beginning of the proof, we have $\alpha_{(n-1),n}=\alpha_{(n-2),n}...=0$. So {\it is the maximum } term is the only one and it corresponds to the maximum in the lexicographic sense of the vector of exponents $\alpha$.

		The fact that $f\in V_{\mu_0}\oplus V_{\mu}\oplus...$ means the following. When using the positive root operators $F_{p,q}$ in some other order, at one of the steps a function is obtained, which is the sum of a highest vector from $V_{\mu}$ and some other summand. This situation occurs exactly when there are terms in the expression \eqref{fpovysh} that are not {\it maximal}.
		
		The weight $weight(f_{\alpha})$ of the resulting $g_{n-k}$-highest vector is calculated as follows. Let
		$h.weight$ denotes $g_n$-the highest weight of the representation in question. Then

		\begin{equation}
			weight(f_{\alpha})=h.weight-\sum_{r<s< -k}\alpha_{r,s}(e_{r}-e_{s}),
		\end{equation}
		
		where $e_{r},e_{s}$ are unit vectors for the corresponding weight components. At the same time, the resulting vector has $weight(f_{\alpha})$ the first $(n-k)$ coordinates are taken.

		When one passes from a maximum term to a non-maximal one to the weight $weight(f_{\alpha})$ a vector of the form $[0,...,1,...,-1,...,0]$  is added.

		The same change in $weight(f_{\alpha})$ occurs when calculating the weight of $weight(f_{\alpha})$ corresponding to the {\it maximum } term, but for the vector of exponents $\delta+r_{\alpha}$, where $r_{\alpha}=e_{-n,...,-r-1,-k}-e_{-n,...,-r-1,-r}-e_{-n,...,-r-1,-s,-k}+e_{-n,...,-r-1,-r,-s}$.
		
		It follows that $a^{\delta}\in \oplus_{s\in\mathbb{Z}^{{\mathcal{K}}}_{\geq 0}}V^{\delta+sr}$.
		
	\end{proof}

	Let us return to the proof of the Theorem \ref{tgc} . Using  Proposition \ref{pfv}, one gets that $\mathcal{G}_{\gamma}(a)$ also belongs to $\oplus_{s\in\mathbb{Z}^{\mathcal{K}}_{\geq 0}}V^{\gamma+sr}$. Since $\mathcal{G}_{\gamma}(a)$ are obtained by the orthogonalization procedure, then $\mathcal{G}_{\gamma}(a)$ is orthogonal to all vectors of the form $\mathcal{F}_{\gamma+sr}(a)$, $s\in\mathbb{Z}^{\mathcal{K}}_{\geq 0}$, $s\neq 0$. Due to the Lemma \ref{lmff}, $\mathcal{G}_{\gamma}(a)$ is also orthogonal to all other $\mathcal{F}_{\delta}(a)$, $\delta\neq\gamma+sr$ $mod\mathcal{B}_{GC}^{g_n}$. So $\mathcal{G}_{\gamma}(a)$ is orthogonal to $\oplus_{s\in\mathbb{Z}^{\mathcal{K}}_{\geq 0},s\neq 0}V^{\gamma+sr}$. It follows that $\mathcal{G}_{\gamma}(a)\in V^{\gamma}$.

	This means that the basis $\mathcal{G}_{\gamma}(a)$ is consistent with the orthogonal decomposition of $V$ into the direct sum of the eigenspaces of $GT$. So $\mathcal{G}_{\gamma}(a)$ is the Gelfand-Tsetlin basis.

	\end{proof}

	Applying now  the same reasoning as in \cite{a4} for the $A$ series, one gets the following statement.

	\begin{theorem}

	The lower-triangular orthogonalization of the $F_{\gamma}$ basis with respect to the order \eqref{poryadok} is a Gelfand-Tsetlin type basis $\mathcal{G}_{\gamma}$.
	
	\end{theorem}

	To prove this statement, we consider the function $G_{\gamma}$ in the variables $A_X$ or $B_X$, representing the vector $\mathcal{G}_{\gamma}$ in A-GKZ realization. Then using the same arguments as for the $A$ series, it is shown that there exists an expression of the form

	$$
	G_{\gamma}=\sum_{s\in\mathbb{Z}_{\geq 0}^{\mathcal{K}} }d_{\gamma}^s\cdot F_{\gamma-sr}
	$$
	
	This expression is an orthogonalization  procedure of Gram-Schmidt.
	The Gram matrix for the $F_{\gamma}$ basis can be explicitly written out using the same  arguments as in \cite{a4} for the $A$ series. Using it one  can explicitly write out the transition matrix from $F_{\gamma}$ to $G_{\gamma}$.

\end{fulltext}

\end{document}